\documentclass[12pt,reqno]{amsart}

\usepackage[utf8]{inputenc}
\usepackage[english]{babel}
\usepackage{amsmath,amssymb,amsfonts,amsbsy,amsthm,amscd,latexsym,graphicx}
\usepackage{empheq,textcomp}

\theoremstyle{definition}
\newtheorem{theorem}{Theorem} [section]
\newtheorem{corollary}[theorem]{Corollary}
\newtheorem{lemma}[theorem]{Lemma}
\newtheorem{proposition}[theorem]{Proposition}
\newtheorem{definition}[theorem]{Definition}

\newtheorem{remark}[theorem]{Remark}
\newtheorem{example}[theorem]{Example}
\newtheorem{question}[theorem]{Question}

\numberwithin{equation}{section}

\newcommand{\C}{\mathbb{C}}

\newcommand{\Fc}{{\mathcal{F}}}
\newcommand{\Gc}{{\mathcal{G}}}
\newcommand{\N}{\mathbb{N}}
\newcommand{\R}{\mathbb{R}}

\newcommand{\Z}{\mathbb{Z}}

\newcommand{\Eq}{\, = \,}

\newcommand{\Le}{\, \le \,}

\newcommand{\Ge}{\, \geq \,}

\newcommand{\qeddef}{{\quad $\diamondsuit$}}

\newcommand{\bigabs}[1]{\bigl|\,#1\,\bigr|}
\newcommand{\Bigabs}[1]{\Bigl|\,#1\,\Bigr|}
\newcommand{\biggabs}[1]{\biggl|\,#1\,\biggr|}
\newcommand{\ip}[2]{\langle\,#1,#2\,\rangle}
\newcommand{\bigip}[2]{\bigl\langle \,#1, \, #2 \,\bigr\rangle}
\newcommand{\Bigip}[2]{\Bigl\langle \,#1, \, #2 \,\Bigr\rangle}
\newcommand{\biggip}[2]{\biggl\langle \,#1, \, #2 \,\biggr\rangle}
\newcommand{\norm}[1]{\|\,#1\,\|}
\newcommand{\bignorm}[1]{\bigl\|\,#1\,\bigr\|}
\newcommand{\Bignorm}[1]{\Bigl\|\,#1\,\Bigr\|}
\newcommand{\biggnorm}[1]{\biggl\|\,#1\,\biggr\|}
\newcommand{\bigparen}[1]{\bigl(\,#1\,\bigr)}
\newcommand{\Bigparen}[1]{\Bigl(\,#1\,\Bigr)}
\newcommand{\biggparen}[1]{\biggl(\,#1\,\biggr)}

\newcommand{\set}[1]{\{#1\}}
\newcommand{\bigset}[1]{\bigl\{\,#1\,\bigr\}}

\newcommand{\parenspace}[1]{(\,#1\,)}

\newcommand{\mpq}{M^{p}(\R)}

\newcommand{\mpqone}{M^{p_1}(\R)}

\newcommand{\Gabor}{\mathcal{G}(g,\alpha,\beta)}
\newcommand{\tfshift}{M_{\beta n}T_{\alpha k}}
\newcommand{\inN}{n\,\in\,\N}
\newcommand{\sumli}{\sum_{n=1}^\infty}
\newcommand{\dousum}{\sum_{k,n\,\in\, \Z}}

\newcommand{\coneab}{C^1_{AC}[a,b]}
\newcommand{\Fgchi}{\Fc(g\cdot \chi_{[a,b]})}
\newcommand{\Fgchic}{\Fc(g\cdot \chi_{[0,c]})}
\newcommand{\gdc}{g\cdot \chi_{[0,c]}}
\begin{document}
\title{Gabor frames with atoms in $M^q(\R)$ but not in $M^p(\R)$ for any $1\leq p<q\leq 2$.}
\author{Pu-Ting Yu}
\address{School of Mathematics, Georgia Institute of Technology,
Atlanta, GA 30332, USA}
\email{pyu73@gatech.edu}
\date{August 2024}

\begin{abstract}
This paper consists of two parts. In the first half, we solve the question raised by Heil as to whether the atom of a Gabor frame must be in $M^p(\R)$ for some $1<p<2$. Specifically, 
 for each $0<\alpha \beta \leq 1$ and $1<q\leq 2$ we explicitly construct Gabor frames $\mathcal{G}(g,\alpha,\beta)$ with atoms in $M^q(\R)$ but not in $M^{p}(\R)$ for any $1\leq p<q$. To construct such Gabor frames, we use box functions as the window functions and show that $$f = \sum_{k,n\in \Z}\ip{f}{\tfshift \Fc(\chi_{[0,\alpha]})} \tfshift ( \Fc(\chi_{[0,\alpha]}))$$ holds for $f\in M^{p,q}(\R)$ with unconditional convergence of the series for any $0<\alpha\beta \leq 1$, $1<p<\infty$ and $1\leq q<\infty$.

 In the second half of this paper, we study two questions related to unconditional convergence of Gabor expansions in modulation spaces. Under the assumption that the window functions are chosen from $M^p(\R)$ for some $1\leq p\leq 2,$ we will prove several equivalent statements that the equation $f = \sum_{k,n\in \Z} \ip{f}{\tfshift \gamma} \tfshift g$ can be extended from $L^2(\R)$ to $M^q(\R)$ for all $f\in M^q(\R)$ and all $p\leq q\leq p'$ with unconditional convergence of the series.  Finally, we characterize all Gabor systems $\set{\tfshift g}_{n,k\in \Z}$ in $M^{p,q}(\R)$ for any $1\leq p,q<\infty$ for which $f = \sum\ip{f}{\gamma_{k,n}} \tfshift g$ with unconditional convergence of the series for all $f$ in $M^{p,q}(\R)$ and all alternative duals $\set{\gamma_{k,n}}_{k,n\in \Z}$ of $\set{\tfshift g}_{n,k\in \Z}$. 
 
\end{abstract}
\maketitle

\section{Introduction}
Fix $\alpha, \beta>0$ and $g\in L^2(\R).$ The \emph{Gabor system} associated with $\alpha,\beta$ and $g$, denoted by $\Gc(g,\alpha,\beta)$, is the sequence $\set{M_{\beta n}T_{\alpha k}g}_{k,n\in \Z}$ consisting of discrete set of time-frequency shifts of $g$, where for fixed $x,\xi \in \R$, $T_x$ is the translation operator $(T_xg)(t)=g(t-x)$ and $M_\xi$ is the modulation operator $(M_{\xi}g)(t)=e^{2\pi i \xi t}g(t)$. If there exist some positive constants $A$ and $B$ such that
\begin{equation}
\label{frame_ineq}
A\,\norm{f}_{L^2(\R)}\Le \dousum \bigabs{\bigip{f}{M_{\beta n}T_{\alpha k}g}_{L^2(\R)}}^2 \Le B\,\norm{f}_{L^2(\R)}
\end{equation}
for all $f\in L^2(\R)$, then we say $\Gabor$ is a \emph{Gabor frame} for $L^2(\R)$ with \emph{frame bounds} $A$ and $B$. If $A=B=1$, then we say $\Gc(g,\alpha,\beta)$ is a \emph{Parseval frame}. The function $g$ is usually called the \emph{atom} or \emph{window function} of $\Gc(g,\alpha,\beta)$. A remarkable property following from Equation (\ref{frame_ineq}) is the existence of another Gabor frame $\Gc(\gamma,\alpha,\beta)$ called the \emph{canonical dual frame} associated with $\Gabor$ such that every $f\in L^2(\R)$ can be expressed as \begin{equation}
\label{Gabor_expan}
f = \dousum \ip{f}{\tfshift \gamma} \tfshift g 
\end{equation}
with unconditional convergence of the series in $L^2$-norm. Equation (\ref{Gabor_expan}) is usually called the \emph{reconstruction formula} or the \emph{Gabor expansion} of $f$. The functions $g$ and $\gamma$ involved in Equation (\ref{Gabor_expan}) are referred to as \emph{atoms} or \emph{window functions}. We refer to \cite{Chr16}, \cite{FS98}, \cite{Gro01} and \cite{Hei11} for background knowledge on the applications and the theory of Gabor frames. A typical question regarding the interplay between Hilbert spaces and Banach spaces is: if certain properties prove to be true for a Hilbert space, then to which Banach spaces can we extend these properties? Consequently, it is natural to ask to which Banach spaces we can extend Equation (\ref{Gabor_expan}) while still preserving unconditional convergence of the series when a Gabor frame is given. 

It was discovered in the early 1980s by Feichtinger and developed in joint work with Gr\"{o}chenig that \emph{modulation spaces} are the right Banach spaces to extend Gabor expansions if the window functions are chosen appropriately (\cite{Fei81}, \cite{FG89a} and \cite{FG89b}). Since then, modulation spaces have been recognized as the appropriate function spaces for time-frequency analysis. Moreover, it has been shown in the past decades that there is a rich interplay between modulation spaces and other existing branches of mathematics such as \emph{pseudodifferential operators} (\cite{GH99}, \cite{GH04}),  \emph{partial differentiable equations} (\cite{BGK20}, \cite{BR16}) and \emph{uncertainty principle} (\cite{GG02}, \cite{HP06}).
Fix a Schwartz function $\psi\in S(\R)$. For $1\leq p,q \leq \infty$, the (unweighted) modulation space $M^{p,q}(\R)$ is the set of all tempered distributions $f\in S'(\R)$ for which $$\norm{f}_{M^{p,q}(\R)} \Eq \Bigparen{\int_\R \Bigparen{\int_\R |V_\psi f(x,w)|^p\,dx }^{q/p}\,dw}^{1/q}<\infty,$$
with the usual modifications if $p=\infty$ or $q=\infty.$ Here $V_\psi f(x,w)$ is the \emph{short-time Fourier transform} with \emph{window function} $\psi$ (see Section 2). Although we only consider modulation spaces over $\R$, modulation spaces were originally defined in a more general setting of \emph{locally compact abelian groups} (\cite{Fei83}). It was proved that modulation spaces are the right Banach spaces to extend Equation (\ref{Gabor_expan}) as long as the window functions are chosen from $M^1(\R)$ (see Theorem \ref{conv_Gabor_modu_original} below). That is, we extend Equation (\ref{Gabor_expan}) by restricting the selection of window functions to functions whose short-time Fourier transforms are $L^1$-functions. Given the close relation between Gabor frames and modulation spaces, it is natural to ask to which level of time-frequency decay the short-time Fourier transform of the atom associated with a Gabor frame can possess. It was raised by Heil in \cite{Hei23} as to whether the atom of a Gabor frame must belong to $M^p(\R)$ for some $1\leq p<2$. That is, must the short-time Fourier transform of a Gabor frame atom decay faster than generic $L^2$-functions? Surprisingly, we will prove that there do exist Gabor frames with atoms not in $M^p(\R)$ for any $1\leq p<2$. Even if we assume the atom is in $M^q(\R)$ for some $1<q<2$, the frame inequality does not necessarily endow the short-time Fourier transform of the atom with a time-frequency decay faster than generic $L^q$-functions.

Another question of main interest in this paper is the limitation of selections of window functions motivated by the extension of Equation (\ref{Gabor_expan}). Specifically, if we choose window functions from $M^p(\R)\setminus M^1(\R)$, then can we still extend Equation (\ref{Gabor_expan}) to $M^{p,q}(\R)$ without losing the unconditional convergence of the series? If the answer is yes, then to which $M^{p,q}(\R)$ we can extend Equation (\ref{Gabor_expan})? Moreover, a Gabor system $\set{\tfshift g}_{k,n\in \Z}$ might admit more than one alternative dual (see Section 5). 
Consequently, it is uncertain whether the unconditional convergence of the series in Equation (\ref{Gabor_expan}) will be preserved when we shift from one alternative dual $\set{\gamma_{kn}}_{k,n\in Z}$ to another even if we are able to extend Equation (\ref{Gabor_expan}) to $M^{p,q}(\R)$. This question, about the \emph{classification of alternative duals}, was studied in \cite{HY23} in the setting of general Hilbert spaces. We will study this question in the setting of general Banach spaces and modulation spaces.

This paper is organized as follows. Notation, terminology, and preliminary results will
be presented in Section 2. To construct the desired Gabor frames, we will need another equivalent norm for certain ranges of modulation spaces. In Section 3, we will derive a family of equivalent norms induced by certain family non-smooth compactly supported functions for $M^{p,q}(\R)$, $1<p<\infty$ and $1\leq q<\infty.$ We then employ these equivalent norms in Section 4 and construct for any fixed $1< q\leq 2$ Gabor frames with atoms in $M^q(\R)$ but not in $M^p(\R)$ for any $1\leq p<q$. In the final section, we will  study the extension of Equation (\ref{Gabor_expan}) in the scenario that the window functions are not necessarily in $M^1(\R)$. In particular, we will present several equivalent statements as to when we can extend Equation (\ref{Gabor_expan}) to certain modulation spaces while still preserving the unconditional convergence of the series when the window functions are not chosen from $M^1(\R)$. Finally, we extend the main result in \cite{HY23} to all Banach spaces that does not contain a topologically isomorphic copy of $c_0$ (see Section 2 for definition). We then use this result to characterize all Gabor systems for which the series in its reconstruction formula converges unconditionally for all elements and all alternative duals in $M^{p,q}(\R)$ for all $1\leq p,q<\infty.$  


\section{Preliminaries}
Throughout this paper, we will use $\ip{\cdot}{\cdot}\colon L^2(\R)\times L^2(\R)\rightarrow \C$ to denote the usual inner product associated with $L^2(\R)$. When $X$ is an arbitrary Banach space with dual space $X^*$, we generalize $\ip{\cdot}{\cdot}$ from $L^2(\R)\times L^2(\R)$ to $X\times X^*$ by interpreting $\ip{x}{x^*}=\overline{\ip{x^*}{x}}$ as $x^*(x)$ for $x^*\in X^*$ and $x\in X.$ We will use $\Fc(f)$ to denote the Fourier transform of $f$ when $\Fc(f)$ is well-defined.

\begin{definition} Let $1\leq p,q\leq \infty$.
\begin{enumerate}
\setlength\itemsep{0.8em}
    \item [\textup{(a)}] The \emph{mixed-norm space} $L^{p,q}(\R^2)$ is space of all (Lebesgue) measurable functions on $\R^2$ for which $$\norm{f}_{L^{p,q}(\R^2)} \Eq \Bigparen{\int_\R \Bigparen{\int_\R | f(x,w)|^p\,dx }^{q/p}\,dw}^{1/q}<\infty,$$
    with the usual modifications if $p=\infty$ or $q=\infty.$
     \item [\textup{(b)}] The \emph{discrete mixed-norm space} $\ell^{p,q}(\Z^2)$ consists of sequences of scalars $(c_{kn})_{k,n\in \Z}$ for which 
     $$\bignorm{(c_{kn})_{k,n\in \Z}}_{\ell^{p,q}(\Z^2)} \Eq \biggparen{\sum_{n\in \Z} \,\Bigparen{\sum_{k\in \Z} |c_{kn}|^p}^{q/p}}^{1/q}<\infty,$$
      with the usual modifications if $p=\infty$ or $q=\infty.$ \qeddef
    \end{enumerate}
\end{definition}

In the case $p=q$, we will simply write $L^{p}(\R)$ or $\ell^p(\Z^2)$ instead of $L^{p,p}(\R)$ or $\ell^{p,p}(\R^2)$. We can now formally define the modulation spaces. 
\begin{definition} Fix a Schwartz function $\psi\in S(\R)$. 
\begin{enumerate}\setlength\itemsep{0.8em}
    \item [\textup{(a)}] Let $\psi\in S'(\R)$ be a tempered distribution. The \emph{short-time Fourier transform} of $f$, denoted by $V_{\psi}f$, is the measurable function on $\R^2$ defined by  $$V_{\psi}f(x,w)\Eq \overline{\ip{M_wT_x \psi}{f}} \Eq  \ip{f}{M_wT_x \psi}.$$

     \item [\textup{(b)}] For $1\leq p,q \leq \infty$, the \emph{modulation space} $M^{p,q}(\R)$ is the space consisting of 
     all tempered distributions $f$ for which $\norm{V_\psi f}_{L^{p,q}(\R^2)}$ is finite, i.e.,
     $$M^{p,q}(\R)=\bigset{f \in S'(\R)\,\big|\, \norm{V_\psi f}_{L^{p,q}(\R^2)}\Eq \norm{f}_{M^{p,q}(\R)}<\infty},$$
      with the usual modifications if $p=\infty$ or $q=\infty.$ If $p=q$, then we simply write $M^p(\R)$ instead of $M^{p,p}(\R)$. \qeddef
    \end{enumerate}
\end{definition}
It is not hard to see that the modulation spaces are Banach spaces under the norm $\norm{\cdot}_{M^{p,q}(\R)}.$
Moreover, it is known that the definition of modulation spaces is independent of the choice of window functions when the choices all come from $S(\R).$
That is, different choices of window functions yield equivalent norms for $M^{p,q}(\R).$ Surprisingly, we can even pick the window function from $M^1(\R)$ and still obtain an equivalent norm for $M^{p,q}(\R).$ We summarize some fundamental properties that are required for this paper in the following lemma. We say that two Banach spaces $X$ and $Y$ are \emph{topologically isomorphic} if there exist a bijective bounded linear operator from $X$ onto $Y$. For notational convenience, we will write $A\lesssim_\phi B$ to mean there exists some constant $C$ depending on $\phi$ such that $A\leq CB.$

\begin{lemma} 
\label{tf_lemmas}
Let $1\leq p,q <\infty.$
\begin{enumerate} \setlength\itemsep{0.8em}
 \item [\textup{(a)}] (\cite[Lemma 11.3.3]{Gro01}) Let $\psi$ be a Schwartz function with $\norm{\psi}_{L^2(\R)}=1$. Then for any $\phi \in S(\R)$ and $f \in S'(\R)$ we have $$|V_\phi f(x,w)|\leq (|V_\psi f|\ast |V_\phi \psi|)(x,w)$$
 for all $x,w\in \R.$

    \item [\textup{(b)}] (\cite[Theorem 11.3.5]{Gro01})  If $p=q$, then $M^p(\R)$ is invariant under Fourier transform. That is, there exists some constant $C>0$  such that $\norm{f}_{M^p(\R)}\lesssim_p \norm{\Fc(f)}_{M^{p}(\R)}$ for all $f\in M^p(\R).$ 
    \item [\textup{(c)}] (\cite[Theorem 11.3.6]{Gro01}) Let $p'$ and $q'$ be the conjugate exponents of $p$ and $q$. Then $(M^{p,q}(\R))^*$ is topologically isomorphic to $ M^{p',q'}(\R).$

    \item [\textup{(d)}] (\cite[Theorem 12.2.2]{Gro01})  For $1\leq p\leq p_1\leq \infty$ and $1\leq q\leq q_1\leq \infty$, we have $$M^{1}(\R)\subseteq M^{p,q}(\R)\subseteq M^{p_1,q_1}\subseteq M^{\infty}(\R)$$
    with $\norm{\cdot}_{M^{p_1,q_1}(\R)}\lesssim_{p,q,p_1,q_1} \norm{\cdot}_{M^{p,q}(\R)}$. 
Moreover, we have that $M^2(\R)=L^2(\R).$
 \item [\textup{(e)}] (\cite[Theorem 11.3.2]{Gro01}) Fix $\psi\in S(\R)$ with $\norm{\psi}_{L^2(\R)}=1.$ Then for any $1\leq p,q<\infty$ we have $ \ip{f}{g} = \ip{V_\psi f}{V_\psi g}$ for any $f\in M^{p,q}(\R)$ and $g\in M^{p',q'}(\R).$
     \item [\textup{(f)}] (\cite[Theorem 1]{FGW92}) For $1\leq p,q<\infty$, $M^{p,q}(\R)$ is topologically isomorphic to $\ell^{p,q}(\Z^2).$
    \end{enumerate}
    \qeddef
\end{lemma}

The following equivalent collection of seminorms for $S(\R)$ will be useful when we are verifying whether a limit of a sequence of tempered distributions is still a tempered distribution. \begin{lemma} (\cite[Corollary 11.2.6]{Gro01})
\label{equiv_norm_Schwartz}
    Fix $g\in S(\R^d)$. The collection of semi-norms $$\norm{V_g f}_{L^\infty_s}\Eq \sup_{x,w\in\R} (1+|x|+|w|)^s|V_g f(x,w)|,\quad s\geq 0, $$ 
    forms an equivalent collection of seminorms for $S(\R).$ 
    \qeddef 
\end{lemma} 

Every Gabor system $\Gabor$ has two associated (not necessarily bounded) operators, the \emph{analysis operator} and the \emph{synthesis operator}. For $1\leq p,q\leq \infty$ the analysis operator denoted by $C_{g}$ is the operator defined by $$C_g(f)\Eq\bigparen{\ip{f}{M_{\beta n}T_{\alpha k} g}}_{k,n\in \Z}$$ on $M^{p,q}(\R)$. The synthesis operator is initially defined on spaces of finite sequences of scalars by $$R_g \bigparen{(c_{kn})_{k,n\in \Z}} \Eq \sum_{k,n\in \Z} c_{kn} M_{\beta n}T_{\alpha k} g.$$ Using the common density argument to extend $R_\gamma$ to $\ell^2(\Z^2)$, Equation (\ref{Gabor_expan}) can be interpreted as the statement that $R_\gamma \circ C_g = I$ on $L^2(\R)$, where $I$ denotes the identity operator.

As indicated in the following theorem, we can now extend Equation (\ref{Gabor_expan}) to certain modulation spaces $M^{p,q}(\R)$ by choosing the window functions from the ``proper" space for which their corresponding analysis operators and synthesis operators are bounded in $M^{p,q}(\R)$ and $\ell^{p,q}(\R)$, respectively. For two arbitrary Banach spaces, $X,Y$, we will write $\norm{\cdot}_{X}\approx \norm{\cdot}_{Y}$ to mean that there exist two positive constants $A,B$ such that $A\norm{\cdot}_{X}\leq \norm{\cdot}_Y \leq B\norm{\cdot}_X$, i.e., $\norm{\cdot}_{X}$ and $\norm{\cdot}_{Y}$ are equivalent norms. Also, recall that we say that a series $\sumli x_n$ \emph{converges unconditionally} in a Banach space $X$ if $\sumli x_{\sigma(n)}$ converges in $X$ for every permutation $\sigma\colon\N\rightarrow\N.$

\begin{theorem} (\cite[pp.\ 258--260]{Gro01})
	\label{conv_Gabor_modu_original}
	Fix $\alpha,\beta >0$. Assume that $g,\gamma \in M^{1}(\R)$. Then the following statements hold.
 \begin{enumerate}
 \setlength\itemsep{0.8em}
     \item  [\textup{(a)}] The analysis operator $C_g$ associated with $g$ is a bounded linear operator from $M^{p,q}(\R)$ to $\ell^{p,q}(\Z^2)$. Moreover, $\norm{C_g}_{op}\lesssim_g \norm{g}_{M^1(\R)}$.
 \item  [\textup{(b)}]  Assume that $f\Eq\sum_{k,n\in \Z} \ip{f}{M_{\beta n}T_{\alpha k}\gamma}  M_{\beta n}T_{\alpha k}g$ converges in $L^{2}(\R)$ with respect to some ordering of the summation. Then 
	for any $1\leq p,q<\infty$, 
	\begin{equation}
	\label{Ga_exp}
	f\Eq\sum_{n,k\in \Z} \ip{f}{M_{\beta n}T_{\alpha k}\gamma}  M_{\beta n}T_{\alpha k}g 
	\end{equation}
	with unconditional convergence of the series in $M^{p,q}(\R)$ for all $f\in M^{p,q}(\R)$.\\
   Furthermore, $\norm{f}_{M^{p,q}(\R)}\approx \bignorm{(\ip{f}{M_{\beta n}T_{\alpha k}\gamma})_{k,n\in \Z}}_{\ell^{p,q}(\Z^{2})}.$\qeddef
    \end{enumerate}
\end{theorem}
By \cite[Theorem 12.1.9]{Gro01}, Theorem \ref{conv_Gabor_modu_original} (b) is essentially saying that we can extend Equation (\ref{Gabor_expan}) while still preserving the unconditional convergence of the series to certain modulation spaces if the window functions are chosen from the ``smallest" Banach space ($M^1(\R)$) that is invariant under translations and modulations. 

A deep question following from this fact is that if the atom of a Gabor frame $\Gabor$ has a certain level of time-frequency decay, then must the atom of its canonical dual frame $\Gc(\widetilde{g},\alpha,\beta)$ have the same level of time-frequency decay? Specifically, if $g\in M^p(\R)$ for some $1\leq p<2$, then must $\widetilde{g}\in M^p(\R)$? The case $p=1$, which we will state below, was confirmed to be true by Gr\"{o}chenig and Leinart in \cite{GL03}, while the case $1<p<2$ is still open. Hereinaftrer, when a Gabor frame $\Gabor$ is given, we let $\widetilde{g}$ denote the atom of the corresponding canonical dual frame.

\begin{theorem} (\cite[Theorem 4.2]{GL03})
	\label{invert_operator_m1}
 Assume that $0<\alpha \beta \leq 1$ and $g\in M^{1}(\R)$ are such that $\Gc(g,\alpha,\beta)$ is a Gabor frame for $L^2(\R)$. Then we have that $\widetilde{g}\in M^1(\R).$

 Consequently, for every $f\in M^{p,q}(\R)$ there exist a sequence of scalars $(c_{kn})_{k,n\in \Z}\in \ell^{p,q}(\Z^2)$ such that $f=\sum_{k,n\in \Z} c_{kn}M_{\beta n}T_{\alpha_k} g$ with unconditional convergence of the series in $M^{p,q}(\R)$ for all $1\leq p,q<\infty.$ In particular, we have that $\norm{(c_{kn})_{k,n\in \Z}}_{\ell^{p,q}(\Z^{2})}\approx \norm{f}_{M^{p,q}(\R)}.$\qeddef
	\end{theorem}

We will use the following technique, named \emph{Painless Nonorthogonal Expansions} due to Daubechies, Grossman, and Meyer, to construct Gabor frames \cite{DGM86}. 

\begin{theorem} (\cite{DGM86})
	\label{painless_non_expan}
	Fix $0< \alpha\beta\leq 1$ and let $g\in L^2(\R^d).$ Assume that supp$(g)\subseteq [0,\beta^{-1}]$. Then $\Gc(g,\alpha, \beta)$ is a frame for $L^2(\R)$ if and only if there exist constants $A,B>0$ such that $$A\beta \Le \sum_{k\in \Z} |g(x-\alpha k)|^2 \Le B\beta\,\,a.e.$$
In this case, $A,B$ are frame bounds for $\Gc(g,\alpha, \beta)$ and $\widetilde{g}=\dfrac{\beta g}{\sum_{k\in \Z} |g(x-\alpha k)|^2}.$\qeddef
\end{theorem}

The atom $g$ in Theorem \ref{painless_non_expan} can be constructed to be smooth as we like, as indicated in the corollary below. 
\begin{lemma} (\cite{DGM86}) 
	\label{window_for_expansion}
	Let $\alpha, \beta >0$ be such that $\alpha\beta < 1.$ Then for any $\alpha <c\leq \frac{1}{\beta}$ there exists some $\psi \in C^{\infty}_c (\R)$ supported in $[0,c]$ such that $\Gc(\psi, \alpha, \beta)$ is a Gabor frame for $L^2(\R).$
\end{lemma}

Finally, we will need the boundedness of the discrete Hilbert transform (for a proof, see \cite[Chapter 13]{FK09}).
\begin{lemma}
	\label{discrete_hilbert}
	For each $1<p<\infty$ there exists a constant $C_p>0$ such that for any sequence of scalars $(c_n)_{\inN} \in \ell^p(\Z)$ we have $$\Bignorm{\Bigparen{\sum_{\substack{n\,\in \,\Z\\ n\neq m}}c_{n} \cdot \frac{1}{m-n}}_{m\in\Z}}_{\ell^p(\Z)}\Le C_p \,\bignorm{(c_n)_{n\in \Z}}_{{\ell^p}(\Z)}.$$
\end{lemma}

\section{Equivalent Norms for $M^{p,q}(\R)$, $1< p<\infty$ and $1\leq q<\infty$}

In this section, we will establish a family of equivalent norms for $M^{p,q}(\R)$ using Fourier coefficients. These equivalent norms will be used frequently in Section 4 when we need to verify whether a function is in certain modulation spaces. For any $a<b$ we let $C^1_{AC}[a,b]$ denote the set of all complex-valued function defined on $\R$ whose restriction to $[a,b]$ has an absolutely continuous derivative on $[a,b].$ That is, 
$$C^1_{AC}[a,b] = \bigset{g\colon \R\rightarrow \C\,\big|\,(g\cdot \chi_{[a,b]})(x) \text{ is differentiable and absolutely continuous on }[a,b]. }.$$

In general, a function in $\coneab$ is not necessary in $M^1(\R)$, and vice versa. For example, for any $a<b$ the box function $\chi_{[a,b]}$ is in $\coneab$ but not in $M^1(\R)$. However, we prove in the next lemma that $\Fc(\coneab \cdot \chi_{[a,b]})\subseteq M^{p,1}(\R)$ for any $p>1$. Here we define $\coneab \cdot \chi_{[a,b]}= \bigset{g\cdot \chi_{[a,b]}\,\big|\,g\in \coneab}.$ 
\begin{lemma}
\label{bdd_of_gchi}
For any $a<b$ we have that $\Fc(\coneab\cdot \chi_{[a,b]})\subseteq M^{p,1}(\R)$ for any $p>1$.
	\end{lemma}
\begin{proof} Since $M^{p,q}(\R)$ is invariant under translation, it suffices to consider $[a,b]=[0,c]$, where $c =b-a.$ By Lemma \ref{window_for_expansion}, there exist some $\beta >0$ and $\psi\in C^{\infty}_c(\R)\subseteq M^1(\R)$ supported in $[0,c]$ for which $\Gc(\psi,\frac{c}{2}, \beta)$ is a frame for $L^2(\R).$ Consequently, $\Gc(\Fc(\psi),\beta, \frac{c}{2})$ is a frame for $L^2(\R).$ Note that since $\bigabs{\ip{M_{\frac{c n}{2}}T_{\beta k}\Fc(\psi)}{\Fc(g\cdot \chi_{[0,c]})}}= \bigabs{\ip{M_{-\beta k}T_{\frac{c n}{2}}\psi}{g\cdot \chi_{[0,c]}}}$, we see that
$\bigabs{\ip{M_{\frac{c n}{2}}T_{\beta k}\Fc(\psi)}{\Fc(g\cdot \chi_{[0,c]})}}$ is nonzero only when $n= -1,0,1.$ 
	For $n=0$, using integration by parts, we obtain
	\begin{align*}
	\begin{split}
	\bigabs{\bigip{T_{\beta k}\Fc(\psi)}{\Fc(g\cdot \chi_{[0,c]})}} \Eq \bigabs{\ip{M_{-\beta k}\psi}{g\chi_{[0,c]}}}
	 &\Eq \biggabs{\int^c_0 e^{-2\pi i \beta k x} \psi(x-c n) \,dx}
	 \Le C \frac{1}{|k|},
	\end{split}
	\end{align*}
for some constant $C>0$ depending on $\phi, \beta.$ The other two cases $n=1$ and $n= (-1)$ follow similarly. Since $\bigparen{ \ip{M_{\frac{\alpha n}{2}}T_{\beta k}(\Fc(\psi)))}{\Fgchi}
 }_{k,n\in \Z} \in \ell^{p,1}(\Z)$ for any $1<p\leq \infty,$ the result then follows by Theorem \ref{conv_Gabor_modu_original} (b). 
	\end{proof}

We will show that every $g\in \coneab$ for which $|g|\geq \delta>0$ on $[a,b]$ for some $\delta>0$ induces a collection of equivalent norms for $M^{p,q}(\R)$ for any $1<p<\infty$ and $1\leq q<\infty.$ The main tools are the following three lemmas. 
\begin{lemma}
	\label{bdd_ana}
	Fix any $\alpha, \beta>0$ and $b>a$.
 Assume that $g\in C^1_{AC} ([a,b])$.  Then  the analysis operator associated with $\Gc(\Fgchi,\beta,\alpha)$, defined by $$C_{\Fgchi}(f)\Eq \Bigparen{\bigip{f}{M_{\alpha n}T_{\beta k}\,\Fgchi}}_{k,n\,\in\, \mathbb{Z}}\,,$$
	is a  bounded linear operator from $M^{p,q}(\R)$ to $\ell^{p,q}(\Z^2)$ for all $1<p<\infty$ and $1\leq q<\infty$.
\end{lemma} 
\begin{proof}
    As before, it suffices to consider $[a,b]=[0,c]$, where $c = b-a.$ Fix $N\in \N$ large enough so that $\frac{\alpha\beta}{N}<\frac{1}{2}$ and $\frac{2\alpha}{N}<c$. By Lemma \ref{window_for_expansion}, there exists a Gabor frame $\Gc(\phi, \frac{\alpha}{N}, \beta)$ for $L^2(\R)$ where the atom $\phi$ is in $C^\infty_c(\R)$ and is supported in $[0,\frac{2\alpha}{N}]$. Consequently, $\Gc(\Fc(\phi), \beta,\frac{\alpha}{N})$ is a Gabor frame for $L^2(\R)$. Since $\Fc(\phi)\in M^1(\R),$
by Theorem \ref{invert_operator_m1}, for every $f\in M^{p,q}(\R)$ there exists a sequence of scalars $(c_{mj})_{m,j\in \Z}\in \ell^{p,q}(\Z^{2})$ for any $f\in M^{p,q}(\R)$  such that $ f = \sum_{m,j\in \Z} c_{mj}M_{\frac{-\alpha j}{N}} T_{-\beta m} (\Fc(\phi)).$ Fix $k,n\in \Z$, we see that 
	\begin{align*}
	\begin{split}
	\Bigabs{\bigip{f}{M_{-\alpha n}T_{-\beta k}(\Fgchic)} }
	&= \Bigabs{\Bigip{\sum_{m,j\in \Z} c_{mj}M_{\frac{-\alpha j}{N}} T_{-\beta m}(\Fc(\phi))}{M_{-\alpha n}T_{-\beta k}\bigparen{\Fgchic}}}\\
	&= \Bigabs{\sum_{m,j\in \Z} c_{mj}\,\Bigip{M_{\frac{-\alpha j}{N}} T_{-\beta m} (\Fc(\phi))}{M_{-\alpha n}T_{-\beta k}\bigparen{\Fgchic}}}\\
	&= \Bigabs{\sum_{m,j\in \Z}\widetilde{c}_{mj}\,\bigip{ M_{-\beta m } T_{\frac{\alpha j}{N}} \phi}{M_{-\beta k}T_{\alpha n}(\gdc)}},
	\end{split}
	\end{align*} 
 where $(\widetilde{c}_{mj})$ is some sequence for which $\norm{(\widetilde{c}_{mj})_{m,j\in \Z}}_{\ell^{p,q}(\Z^{2})}=\norm{c_{mj}}_{\ell^{p,q}(\Z^{2})}$
	We note that since $\Fgchic\in M^{p,1}(\R)$ for any $p>1$ by Lemma \ref{bdd_of_gchi}, we have that the sequence of scalars $$\Bigparen{\bigip{M_{\frac{-\alpha j}{N}} T_{-\beta m}(\Fc(\phi))}{M_{-\alpha n}T_{-\beta k}\bigparen{\Fgchic}}}_{m,j\in \Z} \in \ell^{p,1}(\Z^2)$$ 
 for all $p>1$ by Theorem \ref{conv_Gabor_modu_original} (b).
 So, the interchange of inner product and summation in the second equality above can be justified by using H\"{o}lder's inequality.

 Next, let $$A_{kn} = \Bigabs{\sum_{m\neq k,\,m, j\in \Z}\widetilde{c}_{mj}\,\bigip{ M_{-\beta m } T_{\frac{\alpha j}{N}} \phi}{M_{-\beta k}T_{\alpha n}(\gdc)}}$$ and let $$B_{kn}=\Bigabs{\sum_{j\in \Z}\,\widetilde{c}_{kj}\,\bigip{ M_{-\beta k } T_{\frac{\alpha j}{N}} \phi}{M_{-\beta k}T_{\alpha n}(\gdc)}}$$
	By the Triangle inequality, we see that $\bigabs{\bigip{f}{M_{\alpha n}T_{\beta k}\bigparen{\Fgchic}} } \Le A_{kn}+B_{kn}.$ Let $L\in \N$ be the smallest integer such that $\frac{L\alpha}{N}\geq c.$ Then we see that $\bigip{ M_{-\beta m } T_{\frac{\alpha j}{N}} \phi}{M_{-\beta k}T_{\alpha n}(\gdc)}$ is nonzero only when  $nN-1\leq j\leq nN+L-2$ since $\phi$ is supported in $[0,\frac{2\alpha}{N}]$. Thus, we have $$B_{kn}\lesssim_{\,\phi,g,c}~  \Bigparen{\sum_{nN-1\Le j\Le nN+L-2}|\widetilde{c}_{k,j}|}$$
	Consequently, $\bignorm{(B_{kn})_{k,n\in \Z}}_{\ell^{p,q}(\Z^2)}\Le L\,\norm{f}_{M^{p,q}(\R)}$ for some $L$ depending on $\phi, g$ and $c.$ 

	On the other hand, fix $m,j\in \Z$ and let $D^{k,n}_{m,j} = \bigip{ M_{-\beta m } T_{\frac{\alpha j}{N}} \phi}{M_{-\beta k}T_{\alpha n}(g\chi_{[0,\gamma]})}$. Then by the Triangle Inequality, we obtain $$|A_{kn}|\Le \,\sum_{nN-1\leq j\leq nN+L-2}~\sum_{ m\neq k,\,m\in \Z}  \bigabs{\widetilde{c}_{m,j} D^{k,n}_{m,j}}\,.$$
	Since $m\neq k$, for the case $j=nN$ we compute
	\begin{align*}
	\begin{split}
	D^{k,n}_{m,nN} &\Eq \int^{c+n\alpha}_{n\alpha} e^{ 2\pi i \beta\pi (k-m)x}  g(x-\alpha n)\phi(x-\alpha n) \,dx\\  &\Eq \int^{\frac{2\alpha}{N} }_{0} e^{ 2\pi i \beta (k-m)(x+\alpha n)}  g(x)\phi(x) \,dx \\
	& \Eq \xi_{kn}\nu_{mn} \int^{\frac{2\alpha}{N}}_{0} e^{ 2\pi i \beta (k-m)x} g(x)\phi (x)\, dx,
	\end{split}
	\end{align*}
	where $\xi_{kn} = e^{ 2\pi i k\alpha n}$ and $\nu_{mn} = e^{ -2\pi i \beta m\alpha n}$.
	Using integration by parts, we obtain
	\begin{align*}
	\begin{split}
	\int^{\frac{2\alpha}{N}}_{0} e^{ 2\pi i\beta (k-m)x} (\phi g) (x) dx \lesssim_{\,\beta} \frac{1}{(k-m)}\Bigparen{e^{ 2\pi i \beta (k-m)x} (\phi g) (x)\Big|^{\frac{2\alpha}{N}}_{0}- \int^{\frac{2\alpha}{N}}_{0} e^{2\pi i \beta (k-m)x} (\phi g)' (x) dx}.
	\end{split}
	\end{align*}
	Let $E_{km}= \frac{1}{(k-m)}\Bigparen{e^{ 2\pi i \beta (k-m)x} (\phi g)(x)\big|^{\frac{2\alpha}{N}}_{0}}$ and let $F_{km}= \frac{1}{(k-m)}\int^{\frac{2\alpha}{N}}_{0} e^{ 2\pi i\beta (k-m)x} (\phi g)' (x)\, dx.$ By Lemma \ref{discrete_hilbert}, we see that  
	\begin{align*}
	\begin{split} \bignorm{\bigparen{\sum_{m\neq n,\,m\in\Z}\widetilde{c}_{m,nN} \xi_{kn}\nu_{mn} E_{km}}_{k\in \Z}}_{\ell^p(\Z)}  
	\lesssim_{\phi,g,N,\alpha,\beta} \norm{c_{\cdot,n}}_{\ell^p(\Z)}.
	\end{split}
	\end{align*}
	Using the absolute continuity of $g'$ and integration by parts again, we reach the estimate $$\Bigabs{\frac{1}{ (k-m)} \int^{\frac{2\alpha}{N}}_{0} e^{i2\beta\pi  (k-m)x} (\phi g)' (x) dx}\lesssim_{\phi,g,N,\alpha,\beta} \frac{1}{(k-m)^2}.$$
By Young's convolution inequality, we obtain
	\begin{align*}
	\begin{split} \bignorm{\bigparen{\sum_{m\neq n,\,m\in\Z}\widetilde{c}_{m,nN} \xi_{kn}\nu_{mn} F_{km}}_{k\in \Z}}_{\ell^p}  
	\lesssim_{\phi,g,N,\alpha,\beta}  \norm{c_{\cdot,n}}_{\ell^p(\Z)}.
	\end{split}
	\end{align*}
	Consequently, $\Bignorm{\bigparen{\sum_{m\neq n,\,m\in\Z}\widetilde{c}_{m,nN} D^{k,n}_{m,nN}}_{k\in\Z}}_{\ell^p(\Z)}\lesssim_{\phi,g,N,\alpha,\beta}  \norm{c_{\cdot,n}}_{\ell^p(\Z)}. $ Arguing Similarly for other cases, we see that for any $-1\leq j\leq L+2$, 
	$$\Bignorm{\bigparen{\sum_{m\neq n,\,m\in\Z}\widetilde{c}_{m,nN+\ell} D^{k,n}_{m,nN+j}}_{k\in\Z}}_{\ell^p(\Z)}\lesssim_{\phi,g,N,\alpha,\beta} \norm{c_{\cdot,n}}_{\ell^p(\Z)}.$$
 
	Thus, $$\bignorm{(A_{k,n})_{k,n\in \Z}}_{\ell^{p,q}}\lesssim_{\phi,g,N,\alpha,\beta} \bignorm{\norm{c_{\cdot,n}}_{\ell^p(\Z)}}_{\ell^q(\Z)}\lesssim_{\phi,g,N,\alpha,\beta} \norm{f}_{M^{p,q}(\R)},$$
	and hence, 
	$$\Bignorm{\Bigparen{\bigip{f}{M_{\alpha n}T_{\beta k}(\Fgchic)}}_{k,n\in \Z}}_{\ell^{p,q}(\Z^2)}  \lesssim_{\phi,g,N,\alpha,\beta}\norm{f}_{M^{p,q}(\R)}.$$
\end{proof}
 
\begin{lemma}
	\label{bdd_syn}
	Fix any $\alpha, \beta>0$ and $b>a$.
 Assume that $g\in C^1_{AC} ([a,b])$.  Then  the synthesis operator associated with $\Gc(\Fgchi,\beta,\alpha)$, defined by $$R_{\Fgchi}\bigparen{(c_{kn})_{k,n\in \Z}}\Eq \sum_{k,n\in \Z} c_{kn} M_{\alpha n}T_{\beta k}\bigparen{\Fgchi},$$
	is a bounded linear operator from $\ell^{p,q}(\Z^2)$ to $M^{p,q}(\R)$ for all $1<p<\infty$ and $1\leq q<\infty$. 
 
 Consequently, $\sum_{k,n\in \Z} c_{kn} M_{\alpha n}T_{\beta k}(\Fgchic)$ converges unconditionally in $M^{p,q}(\R)$ for all $(c_{kn})_{k,n\in \Z}\in \ell^{p,q}(\Z^2)$.
\end{lemma} 
\begin{proof}
    Similar to Theorem \ref{bdd_ana}, it suffices to consider $[a,b]=[0,c]$, where $c = b-a.$
 	We first  show that $R_{\Fgchic}\bigparen{(c_{kn})_{k,n\in \Z}}\in S'(\R)$ for any $(c_{kn})_{k,n\in \Z}\in \ell^{p,q}(\Z^2)$. Let $\phi \in S(\R)$ be an arbitrary Schwartz function and let $ r\geq 0$ be large enough that $\bignorm{(1+|x|+|w|)^{-r}}_{L^1(\R^2)}<\infty.$ Fix $\varphi\in S(\R)$. Note that $V_\varphi \phi\in S(\R)$ by \cite[Theorem 11.2.5]{Gro01}. Then by Theorem \ref{conv_Gabor_modu_original}(b) and Lemma \ref{bdd_of_gchi} we compute 
  \begin{align*}
      \begin{split}
          \Bigabs{\Bigip{\phi}{\sum_{k,n\,\in\, \mathbb{Z}} c_{kn} M_{\alpha n} T_{\beta_{k}} \bigparen{\Fgchic}}}&\Le \sum_{k,n\,\in\, \mathbb{Z}} \bigabs{c_{kn}}\,\Bigabs{\Bigip{\phi}{  M_{\alpha n} T_{\beta_{k}} \bigparen{\Fgchic}}}\\
          &\Le \norm{c_{kn}}_{\ell^{p,q}(\Z^2)} \Bignorm{\Bigip{\phi}{  M_{\alpha n} T_{\beta_{k}}\bigparen{\Fgchic}}}_{\ell^{p',q'}(\Z^2)}\\
          &\lesssim_{(c_{kn})} \norm{\phi}_{M^1(\R)} \norm{\Fgchic}_{M^{p',q'}(\R)}\\
          &\lesssim_{(c_{kn}),\,g}  \bignorm{\underbrace{V_\varphi \phi}_{\in S(\R)}}_{L^1(\R^2)}\\
          & \lesssim_{(c_{kn}),\,g\,} \bignorm{ |V_\gamma\phi|(1+|x|+|w|)^r(1+|x|+|w|)^{-r}}_{L^1(\R^2)}\\
          &\lesssim_{(c_{kn}),\,g,\,r} \sup_{x,w\,\in \R} |V_\varphi \phi|(1+|x|+|w|)^r.
      \end{split}
  \end{align*}
Consequently, $\sum_{k,n\,\in\, \mathbb{Z}} c_{kn} M_{\alpha n} T_{\beta k} (\Fgchic)\in S'(\R)$ by Lemma \ref{equiv_norm_Schwartz}. The rest of the proof follows from a similar idea to the proof of Lemma $\ref{bdd_ana}$. Fix $N$ large enough so that $\frac{\alpha\beta}{N}<\frac{1}{2}$ and $\frac{2\alpha}{N}<c$. By Lemma \ref{window_for_expansion}, there exists $\psi\in C^\infty_c(\R)\subseteq M^1(\R)$ supported in $[0,\frac{2\alpha}{N}]$ for which $\Gc(\Fc(\psi), \beta,\frac{\alpha}{N})$ is a Gabor frame for $L^2(\R)$. By Theorem \ref{conv_Gabor_modu_original}, it suffices to show that 
$$\Bignorm{C_{\Fc(\psi)} \Bigparen{\sum_{k,n\,\in\, \mathbb{Z}} c_{kn} M_{-\alpha n} T_{-\beta k} \bigparen{\Fgchic}}}_{M^{p,q}(\R)}\Le C\,\norm{(c_{kn})_{k,n\in \Z}}_{\ell^{p,q}(\Z^2)},$$
for some constant $C>0$, where $C_{\Fc(\psi)}$ is the analysis operator associated with $\Gc(\Fc(\psi), \beta,\frac{\alpha}{N})$.
 Using Lemma \ref{bdd_of_gchi}, Theorem \ref{conv_Gabor_modu_original}, and H\"{o}lder's inequality to justify the interchange of inner product and summation in $C_{\Fc(\psi)} \bigparen{\sum_{k,n\,\in\, \mathbb{Z}} c_{kn} M_{-\alpha n} T_{-\beta k} \bigparen{\Fgchic}}$, it remains to show 
  \begin{align*}
 \begin{split}
  \biggabs{\sum_{k,n\in \Z} \widetilde{c}_{kn}\Bigip{M_{-\beta k} T_{\alpha n}(g\cdot\chi_{[0,c]})}{M_{-\beta m}T_{\frac{\alpha j}{N}}\psi}} \Le C\,\norm{(c_{kn})_{k,n\in \Z}}_{\ell^{p,q}(\Z^2)},
	\end{split}
	\end{align*} 
 where $\widetilde{c}_{kn}= c_{kn}e^{i2\pi \alpha \beta nk/N}$. Arguing similarly to the proof used for Lemma \ref{bdd_ana}, we reach the estimate
 \begin{equation}
 \label{bddineq_syn}
     \bignorm{R_{\Fgchi}\bigparen{(c_{kn})_{k,n\in \Z}}}_{M^{p,q}(\R)}\lesssim_{g,a,b,\psi,N} \norm{(c_{kn})_{k,n\in \Z}}_{\ell^{p,q}(\Z^2)}.
\end{equation}
To see that $R_{\Fgchi}\bigparen{(c_{kn})_{k,n\in \Z}}$ converges unconditionally in $M^{p,q}(\R)$, we will use an equivalent condition (\cite[pp.\ 98]{Gro01}) of unconditional convergence. 
For any $\varepsilon>0$ there exists a finite subset $F_0\subseteq \Z^{2}$ such that $$\norm{(c_{kn})_{k,n\in \Z}-(c_{kn})_{k,n\in \Z}\cdot \chi_{F_0}}_{\ell^{p,q}(\Z^{2})}<\varepsilon.$$ Here $(c_{kn})_{k,n\in \Z}-(c_{kn})_{k,n\in \Z}\cdot \chi_{F_0}$ means the sequence $(d_{kn})_{k,n\in \Z}$ for which $d_{kn}=0$ if $(k,n)\in F_0$ and $d_{kn}=c_{kn}$ if $(k,n)\notin F_0$.
	Then for any $F\supseteq F_0$ there is a constant $C'>0$ such that $$\norm{R_{\Fgchi}\bigparen{(c_{kn})_{k,n\in \Z}-(c_{kn})_{k,n\in \Z}\cdot \chi_{F}}}_{M^{p,q}(\R)}\leq C' \varepsilon,$$
 for some $C'>0.$ Thus, $R_{\Fgchi}\bigparen{(c_{kn})_{k,n\in \Z}}$ converges unconditionally in $M^{p,q}(\R)$.
\end{proof}

By Theorem \ref{painless_non_expan}, every $g\in \coneab$ for which $|g|$ is bounded below by a positive number induces a collection of Gabor frames $\Gabor$ within certain ranges of the translation and modulation parameters. We just established the boundedness of the analysis operator associated with $\Gabor$ in $M^{p,q}(\R)$ for all $1<p<\infty$ and $1\leq q<\infty$ in Lemma \ref{bdd_ana}. It remains to show that the synthesis operator associated with the canonical dual frame $\Gc(\widetilde{g},\alpha,\beta)$ of $\Gabor$ is also bounded in $M^{p,q}(\R)$ for all $1<p<\infty$ and $1\leq q<\infty$. Then we can obtain a collection of equivalent norms for $M^{p,q}(\R)$ via the following inequality
\begin{equation} 
\label{operator_ineq}
\norm{f}_{M^{p,q}(\R)} \Eq \norm{(R_{\widetilde{g}}\circ C_g)(f)}_{M^{p,q}(\R)} \lesssim_{\,\widetilde{g}} \norm{C_g(f)}_{\ell^{p,q}(\Z^2)}\lesssim_{\,\widetilde{g},g} \norm{f}_{M^{p,q}(\R)}.  
\end{equation}
The next lemma shows that the atom of the canonical dual frame associated with $\Gabor$ can be written as a sum of finitely many functions in $\coneab.$
\begin{lemma}
\label{decom_frame_operator}
Assume that $g\in C^1_{AC}[a,b]$ and $|g|\geq \delta >0$ on $[a,b].$ Then for any $b-a>\alpha >0$ there exist finitely many non-overlapping closed intervals $I_\ell = [a_\ell,b_\ell]$ such that \begin{enumerate}
    \setlength\itemsep{0.8em}
        \item [\textup{(a)}]
$ [a,b]=\bigcup_{\ell=1}^N I_\ell$
\item [\textup{(b)}] For each $1\leq \ell\leq N$ there exists $h_\ell\in C^1_{AC}(I_\ell)$ for which $h_\ell \Eq \displaystyle \frac{g\cdot \chi_{I_\ell}}{ \sum_{n\in \Z}|T_{\alpha n}(g\cdot \chi_{[a,b]})|^2}$ everywhere on $(a_\ell,b_\ell).$
    \end{enumerate}

\end{lemma}
\begin{proof}
The case $\alpha=b-a$ simply follows by letting $N=1$ and $I_1 = [a,b].$
    If $\alpha < b-a$, then we let $K=\lceil\frac{b-a}{\alpha}\rceil$ and define intervals $$R_1=[a,a+\alpha],R_2=[a+\alpha,a+2\alpha],\dots,R_K=[a+(K-1)\alpha,b],$$
    and 
    $$L_1=[b-\alpha,b],L_2=[b-2\alpha,b-\alpha],\dots,L_K=[a,b-(K-1)\alpha].$$
    For each $1\leq i, j\leq K$, we see that $\sum_{n\in \Z}|T_{\alpha n}(g\cdot \chi_{[a,b]})|^2 \cdot \chi_{R_i\cap L_j}$ only has finitely many nonzero terms.
    Moreover, since $|g|\geq \delta$ on $[a,b]$, we see that the only points of discontinuity of $\frac{g\cdot \chi_{R_i\cap L_j}}{ \sum_{n\in \Z}|T_{\alpha n}(g\cdot \chi_{[a,b]]})|^2}$ are the two endpoints of $R_i\cap L_j$. Using the fact that $|g|\geq \delta$ again, we obtain some function $h_{ij}\in \coneab$ by continuously extending $\frac{g\cdot \chi_{R_i\cap L_j}}{ \sum_{n\in \Z}|T_{\alpha n}(g\cdot \chi_{[a,b]})|^2}$ from the interior of $R_i\cap L_j$ to the endpoints of $R_i\cap L_j$. 
\end{proof}

\begin{theorem}
 		\label{conv_Gabor_ext} Fix $a<b$ and fix $0<\alpha\leq (b-a)\leq \beta^{-1}$. Assume that $g\in \coneab$ for which  $|g|\geq \delta >0$ on $[a,b].$  Then for any $1<p<\infty,\,1\leq q<\infty$ the following statements hold:    
    \begin{enumerate}
       \setlength\itemsep{0.8em}
       \item [\textup{(a)}] $\Gc(g\cdot \chi_{[a,b]},\alpha,\beta)$ is a Gabor frame for $L^2(\R)$. 
    \item [\textup{(b)}] The Gabor expansions \begin{equation}
      \label{Gabor_exp}
    \begin{split}
        f&\Eq\sum_{n,k\in \Z} \bigip{f}{M_{\alpha n}T_{\beta k}(\Fc(\,\widetilde{g}\,)}  \,M_{\alpha n}T_{\beta k} \bigparen{\Fgchi}\\
        &\Eq\sum_{n,k\in \Z} \bigip{f}{M_{\alpha n}T_{\beta k} \bigparen{\Fgchi}}  M_{\alpha n}T_{\beta k} \bigparen{\Fc(\,\widetilde{g}\,)},
    \end{split}
     \end{equation}
     hold with unconditional convergence of the series in $M^{p,q}(\R)$ for all $f\in M^{p,q}(\R)$, where $\widetilde{g}$ denotes the atom of the canonical dual frame associated with $\Gc(g\cdot \chi_{[a,b]},\alpha,\beta)$. 
     \item [\textup{(c)}]
 			    The following norm equivalence holds:  $$\Bignorm{\Bigparen{\bigip{f}{M_{\alpha n}T_{\beta k}(\Fgchi)}}_{k,n\in \Z}}_{\ell^{p,q}(\Z^2)} \approx \norm{f}_{M^{p,q}(\R)}.$$
         \end{enumerate}
\end{theorem}
\begin{proof} For notational convenience, we will write $g\cdot \chi_{[a,b]} = g_{a,b}$. 
    Statement (a) simply follows from Theorem \ref{painless_non_expan}. Note that the atom $\widetilde{g}$ of the canonical dual frame associated with $\Gc(g_{a,b},\alpha,\beta)$ is $\widetilde{g}=\displaystyle \frac{\beta g_{a,b}}{ \sum_{n\in \Z}|T_{\alpha n}(g_{a,b})|^2}$.
    
    For statement (b), the synthesis operator $R_{\Fc(\widetilde{g})}\colon\ell^{p,q}(\Z^2)\rightarrow M^{p,q}(\R)$ and the analysis operator $C_{\Fc(\widetilde{g})}\colon M^{p,q}(\R)\rightarrow \ell^{p,q}(\Z^2)$ associated with $\Gc(\Fc(\widetilde{g}),\beta,\alpha)$ are bounded by Lemma \ref{decom_frame_operator}, Theorem \ref{bdd_ana}, Theorem \ref{bdd_syn}, and the Triangle Inequality. Statement (b) then follows from the fact that $\Gc(\Fc(g_{a,b}),\beta,\alpha)$ is a frame for $L^2(\R)$ and
    $$R_{\Fc(\widetilde{g})}\circ C_{\Fc(g_{a,b})} \Eq I_{L^2(\R)} \Eq  R_{\Fc(g_{a,b})}\circ C_{\Fc(\widetilde{g})}.$$

    The proof of statement (c) now simply follows from Inequality (\ref{operator_ineq}).
\end{proof}
\section{Gabor frames with atoms in $M^{q}(\R)$ but not in $M^{p}(\R)$ for any $1\leq p< q\leq 2$}
We will construct Gabor frames with atoms not in $M^p(\R)$ for any $p<2$ in this section. Even more, for each $1<q\leq 2$ we will construct Gabor frames with atoms in $M^q(\R)$ but not in $M^p(\R)$ for any $1\leq p<q.$ We will need the following corollary, which is a special case of Theorem \ref{conv_Gabor_ext}. We mention that some results that characterize compactly supported elements in $M^{p,q}(\R)$ using Fourier transforms and Fourier coefficients can be found in \cite{BO11},\cite{MK06} (see also \cite{BO20}).

\begin{corollary}
\label{equi_norm_Fcoef}
	For any $1<p\leq 2$ and any $0< \alpha \beta \leq 1$, we have the following norm equivalence:
 $$\norm{f}_{M^{p}(\R)}\approx \Bigparen{\sum_{k,n\in \Z} \bigabs{\Fc(f\cdot\chi_{[\alpha n,\alpha(n+1)]})(\beta k)}^p}^{1/p}$$ for all $f\in M^{p}(\R)$.
 
 Consequently, for any bounded interval $I$ the multiplication operator $T(f)=f \chi_I$ is a bounded linear operator from $M^p(\R)$ to $M^p(\R)$.
	\end{corollary}
\begin{proof}
	By Corollary \ref{conv_Gabor_ext}, we have $\bignorm{\bigparen{\bigip{f}{M_{\alpha n}T_{\beta k}(\Fc(\chi_{[0,\alpha]})}}_{k,n\in \Z}}_{\ell^{p,q}(\Z^2)} \approx \norm{f}_{M^{p,q}(\R)}.$ Since $M^{p}(\R)\subseteq M^2(\R) = L^{2}(\R)$ and $M^p(\R)$ is invariant under Fourier transform, we see that
	
	\begin{align*}
\begin{split}
	\norm{f}_{M^{p}(\R)}\approx \bignorm{\Fc(f)}_{M^{p}(\R)} &\approx  \Bignorm{\Bigparen{\bigip{\Fc(f)}{M_{\alpha n}T_{\beta k}\bigparen{\Fc(\chi_{[0,\alpha]})}}}_{k,n\in \Z}}_{\ell^{p}(\Z^2)}\\ &\approx  \Bignorm{\bigparen{\bigip{f}{M_{\beta k}T_{-\alpha n}\chi_{[0,\alpha]}}}_{k,n\in \Z}}_{\ell^{p}(\Z^2)}.
	\end{split}
	\end{align*}
 Finally, applying translation to $f\cdot \chi_I$ if necessary, we may consider $I=[0,c]$ for some $c>0.$ Choosing $\beta$ small enough that $c < \frac{1}{\beta},$ we then obtain 
 \begin{align*}
 \begin{split}
     \bignorm{f\cdot\chi_{[0,c]}}_{M^p(\R)} &\approx \,\Bignorm{\Bigparen{\bigip{f\cdot \chi_{[0,c]}}{M_{\beta k}T_{c n}\chi_{[0,c]}}}_{k,n\in \Z}}_{\ell^{p}(\Z^2)}\\
     &\leq \,\Bignorm{\Bigparen{\bigip{f}{M_{\beta k}T_{c n}\chi_{[0,c]}}}_{k,n\in \Z}}_{\ell^{p}(\Z^2)}\\
     &\approx \,\norm{f}_{M^{p}(\R)}.
     \end{split}
 \end{align*}
 \end{proof}
\begin{remark}
We remark that for $d>1$ we also have  $$\norm{f}_{M^{p}(\R^d)}\approx \Bigparen{\sum_{k=(k_1,..,k_d),\,n=(n_1,..,n_d)\in \Z^d} \bigabs{\Fc(f\cdot \chi_{I_n})(\beta k_1,...,\beta k_d)}^p}^{1/p}$$
for any $\alpha \beta \leq 1,$
where $I_n=[\alpha n_1,\alpha(n_1+1)]\times\cdots \times [\alpha n_d,\alpha(n_d+1)].$
 Since this is not the main interest of this paper, we provide a sketch of the proof here.  
 First, Theorem \ref{painless_non_expan} remains true in $\R^d$. Consequently, we can show that for any $0< \alpha\beta \leq 1$ and $0<c\leq \beta^{-1}$ there exists some $\psi\in C_c(\R)$ supported in $[0,c]$ such that $\Gc(\underbrace{\psi \otimes\cdots \otimes \psi}_{d\text{-times}},\alpha,\beta)$ is a Gabor frame for $L^2(\R^d)$ by using a similar proof to 
 Lemma \ref{window_for_expansion}. Here $\psi \otimes\cdots \otimes \psi$ denotes the function $\phi$ defined on $\R^d$ for which $\phi(x_1,\dots,x_d)=\psi(x_1)\cdots\psi(x_d).$  Next, by using Lemma \ref{discrete_hilbert} step by step, we obtain that
 $$\bignorm{\sum_{\substack{n_1,...,n_d\,\in \,\Z\\ (n_1,\dots,n_d)\neq (m_1,\dots ,m_d)}}c_{n_1,\dots,n_d} \cdot \frac{1}{m_1-n_1}\cdots\frac{1}{m_d-n_d}}_{\ell^p(\Z^d)}\Le C_p \norm{c}_{{\ell^p}(\Z^d)}$$
 for some $C_p>0$ and any $c=(c_{n_1,\dots,n_d})_{n_1,\dots,n_d\in \Z}\in \ell^p(\Z^d)$ for all $1<p<\infty.$ Finally, arguing similarly to Lemma \ref{bdd_ana}, Theorem  
 \ref{bdd_syn} and Theorem \ref{conv_Gabor_ext}, we can obtain the $\R^d$-version of Corollarly \ref{equi_norm_Fcoef}. 

 \begin{remark}
      Corollary \ref{equi_norm_Fcoef} cannot be extended to $M^1(\R)$. A counterexample is simply $\chi_{[0,c]}$ for any $c>0$, which is in $M^{1,p}(\R)\setminus M^{1,1}(\R)$ for any $1<p\leq \infty.$
 \end{remark}
 
\end{remark}

The construction of our examples relies on a variant of the Shapiro-Rudin polynomials. Recall that the Shapiro-Rudin polynomials are defined by the following recursion. Let $P_0=1=Q_0$. For $n\geq 1$, $P_n$ and $Q_n$ are defined by 
\begin{align*}
    P_n &= P_{n-1}+e^{i2\pi2^{n-1}x}Q_{n-1},\\
    Q_n &= P_{n-1}-e^{i2\pi2^{n-1} x}Q_{n-1}.
\end{align*}
Next, for each $n\geq 1$ we define $f_n = (P_{n}-P_{n-1})\cdot \chi_{[0,1]}$. Note that for each $n\geq 1$ the function $f_n$ has the following properties (for example, see \cite{Kat04}):
\begin{enumerate}
\setlength\itemsep{0.8em}
    \item [\textup{(a)}] $f_n = \sum_{m=2^{n-1}}^{2^n-1}\epsilon_m e^{i2\pi mx}, \text{\,where\,} \epsilon_m = 1 \text{ or} -1, $
    \item [\textup{(b)}] $\sumli|\ip{f_j}{M_nT_k\chi_{[0,1]}}|^m = \sumli|\Fc(f_j)(n)|^m = 2^{j-1}$ \quad \text{for all $m>0$},
    \item [\textup{(c)}] $\norm{f_n}_{L^\infty([0,1])}\Le 2^{\frac{n+1}{2}}.$
\end{enumerate}
\begin{lemma}
\label{bdd_fn_not_in_Ap}
 Fix $1<p<q\leq 2$.
    Let $0\leq a<b\leq 1$ be such that $b-a=\frac{1}{L}$ for some $L\in \N$. Then for any $\epsilon>0$  there exists a function $h$ supported in $[a,b]$ for which \begin{enumerate}
    \setlength\itemsep{0.8em}
        \item [\textup{(a)}] $|h|\leq 1$ on $[a,b],$
     \item [\textup{(b)}] $\norm{h}_{M^q(\R)} < \epsilon,$ 
     \item [\textup{(c)}] $\norm{h}_{M^p(\R)} = \infty.$ 
    \end{enumerate}
\end{lemma}
\begin{proof}
    For each $n\geq 0$ we define $g_p(x)= \sum_{n=1}^\infty 2^{-n/p}f_n.$ Since $\norm{f_n}_\infty\leq  2^{\frac{n+1}{2}}$, we see that $\norm{g_p}_{\infty}\leq 2^{1/2}\sum_{n=1}^\infty 2^{n(1/2-1/p)}<\infty.$ By Corollary \ref{equi_norm_Fcoef} we compute 
    \begin{align*}
        \begin{split}
            \norm{g_p}_{M^q(\R)} &\approx \sum_{n,k\in \Z}\,\bigabs{\ip{g_p}{M_{k}T_{n}\chi_{[0,1]}}}^q \\
            &\Eq \sum_{k\in \N} 2^{-\frac{qk}{p}}(2^k-1-2^{k-1}+1) \\
            &\Eq 2^{-1}\sum_{k\in \N} 2^{k(1-\frac{q}{p})}\\
            &<\infty.
        \end{split}
    \end{align*}
    Since $M^q(\R)$ is invariant under translation, the multiplication by $\chi_{[a,b]}$, and dilation $x\mapsto Kx$ for any $K\in \N$, we can find $M$ large enough that $\frac{\norm{g_p}_\infty}{M}<1$ and $\frac{\norm{g_p(L(x-a))\chi_{[a,b]}}_{M^q(\R)}}{M}<\epsilon.$ We then see that $h(x)=\dfrac{g_p(L(x-a))}{M}\chi_{[a,b]}$ satisfies $(a)$ and $(b)$. It remains to show that $\norm{h}_{M^p(\R)}=\infty.$
    For each $n\in \N$ the function $f_n=\sum_{m=2^{n-1}}^{2^n-1}\epsilon_m e^{i2\pi mx}$ contains at least $\lfloor 2^{n-1}/L \rfloor$ exponential terms $e^{i2\pi mx}$ for which $L|m$. So, we have 
\begin{equation*}
\sum_{L|n,n\in\N} \bigabs{\Fc(g_p)(n/L)}^p \Ge \sum_{L|n,\,n\in \N}2^{-n} (\lfloor 2^{n-1}/L \rfloor)\Ge \sum_{L|n,\,n\geq L}2^{-n} 2^{n-1-L}\Eq \infty. 
\end{equation*}
Finally, since 
$$|\Fc(h)(n)|= \frac{1}{M}\Bigabs{\int^{b}_a g_p\bigparen{L(x-a)}e^{-i2\pi nx}\,dx}=\Bigabs{\int^{1/L}_0 g_p(Lx)e^{-i2\pi nx}\,dx} = \frac{M}{L}\bigabs{\Fc(g_p)(n/L)},$$
we obtain that
\begin{align*}
        \begin{split}
            \norm{h}^p_{M^p(\R)} &\approx \displaystyle\sum_{\substack{n\in \Z}}|\Fc(h)(n)|^p\Ge \sum_{\substack{L|n,\\n\in \Z}}|\Fc(h)
(n)|^p \Eq \frac{M}{L}\sum_{\substack{L|n,\\n\in \Z}}\bigabs{\Fc(g_p)(\frac{n}{L})}^p=\infty
        \end{split}
    \end{align*}
    by Corollary \ref{equi_norm_Fcoef}.
\end{proof}
\begin{theorem}
\label{counter_example}
For each $0<\alpha\beta \leq 1$ and $1<q\leq 2$ there exists a function $g\in M^q(\R)$ for which $\Gc(g,\alpha,\beta)$ is a frame for $L^2(\R)$ but $g\notin M^{p}(\R)$ for any $1\leq p<q$.

Moreover, if $q=2$, then for each $0<\alpha\beta \leq \frac{1}{2}$ there exists $g\in M^2(\R)$ for which $\Gc(g,\alpha,\beta)$ is a Parseval frame for $L^2(\R)$ but $g\notin M^p(\R)$ for any $1\leq p<2.$ 
\end{theorem}
\begin{proof}
Using dilation, it is enough to consider the case 
  $\alpha =1$ and $\beta <1$. Let $(p_k)_{k\in \N}$ be a sequence contained in $(1,q)$ that increases to $q$. Let $I_1=[0,1/2]$ and for each $k\geq 1$ let $I_k=[\,\sum_{n=1}^k 2^{-n},\sum_{n=1}^{k+1} 2^{-n}\,]$ (Note that $\cup_{k=1}^\infty I_k=[0,1)$). 
    By Lemma \ref{bdd_fn_not_in_Ap}, for each $k\in \N$ there exists a function $g_k$ supported in $I_k$ which satisfies \begin{enumerate}
    \setlength\itemsep{0.8em}
        \item [\textup{(a)}] $|g_k|\leq 1$ on $I_k$, 
        \item [\textup{(b)}] $\norm{g_k}_{M^q(\R)}\leq \frac{1}{2^k},$
        \item [\textup{(c)}] $\norm{g_k}_{M^{p_k}(\R)}=\infty.$ \end{enumerate}
    
    Let $g=\sum_{k=1}^\infty g_k+2\chi_{[0,1]}$. Then we see that $\text{supp}(g)=[0,1]$ and $1\leq |g|\leq 3$ for almost every $x\in [0,1]$.  By Theorem \ref{painless_non_expan}, $\mathcal{G}(g,1,\alpha\beta)$ is a Gabor frame for $L^2(\R)$. Moreover, by Lemma \ref{bdd_of_gchi}, we have that $$\norm{g}_{M^q(\R)}\Le \sum_{k\in\N} \norm{g_k}_{M^q(\R)}+ \norm{\chi_{[0,1]}}_{M^{q}(\R)}<\infty. $$

    It remains show that $g\notin M^p(\R)$ for any $1\leq p<q.$ Equivalently, we show that $g\notin M^{p_k}(\R)$ for any $k\in \N$ instead. If $g$ was in $M^{p_k}(\R)$ for some $k\in \N$, we would have 
\begin{equation}
\label{tri_ineq1}
    \sum_{k=1}^\infty g_k =g-2\chi_{[0,1]} \in M^{p_k}(\R),
    \end{equation}
    by the Triangle Inequality. But then we must have $(T_{-\sum_{n=1}^k 2^{-n}}) (\sum_{k\in\N} g_k) \in M^{p_k}(\R)$. However, by considering $\alpha = 2^{-(k+1)}$ and $\beta =1$ in Corollary \ref{equi_norm_Fcoef}, we see that 
    \begin{align*}
        \begin{split}
            \bignorm{(T_{-\sum_{n=1}^k 2^{-n}})g}_{M^{p_k}(\R)} &\approx \biggnorm{\Bigparen{\bigip{T_{-\sum_{n=1}^k 2^{-n}})g}{M_ mT_{-2^{-(k+1)}j}\chi_{[0,2^{-(k+1)}]}}}_{m,j\in \Z}}_{\ell^{p_k}(\Z^2)}\\
            &\Ge \biggnorm{\Bigparen{\bigip{T_{-\sum_{n=1}^k 2^{-n}})g}{M_ m\chi_{[0,2^{-(k+1)}]}}}_{m\in\Z}}_{\ell^{p_k}(\Z)}\\
 &\Eq \biggnorm{\Bigparen{\bigip{g_k\cdot \chi_{[0,2^{-(k+1)}]}}{M_ mT_{-2^{-(k+1)}j} \, \chi_{[0,2^{-(k+1)}]}}}_{m\in\Z}}_{\ell^{p_k}(\Z)}\\
            &\approx \bignorm{g_k}_{M^{p_k}(\R)}\\  
            &\Eq \infty,
        \end{split}
    \end{align*}
which is a contradiction. Therefore, $g\notin M^{p}(\R)$ for any $1\leq p<q.$    

If $0<\alpha\beta\leq \frac{1}{2}$, similarly, then it is enough to consider the case $\alpha=1$ and $\beta\leq \frac{1}{2}.$ Let $\set{p_k}_{\inN} \subseteq (1,2)$ be a sequence that increases to $2.$ 
Using a similar construction above, we can find $h$ supported in $[0,1]$ for which 
\begin{enumerate}
    \setlength\itemsep{0.8em}
        \item [\textup{(a)}] $1\leq |h|\leq 3$ in $[0,1]$, 
        \item [\textup{(b)}] $\norm{h\cdot\chi_{[\sum_{n=1}^k 2^{-n},\sum_{n=1}^{k+1} 2^{-n}]}}_{M^{p_k}(\R)}=\infty ~~ \text{for all $k\in \N$}.$ \end{enumerate}
Dividing $h$ by a large constant if necessary, we may assume that $\delta<|h|^2\leq \beta$ in $[0,1]$ for some $\delta >0$. Next, we define
$$g(x)\Eq h(x) + \sqrt{\beta-|h(x)|^2} \cdot \chi_{[1,2]}
$$
Then we see that $g\in L^2(\R)$ and $\sum_{k\in \Z}|g(x-k)|^2 = \beta $ a.e. Consequently, $\Gc(g,1,\beta)$ is a Parseval frame for $L^2(\R)$ by Theorem \ref{painless_non_expan}. To see that $g \notin M^{p}(\R)$ for any $1<p<2$, similarly, by considering $\alpha = 2^{-(k+1)}$ and $\beta =1$ in Corollary \ref{equi_norm_Fcoef}, we see that for each $k\in \N$ 
    \begin{align*}
        \begin{split}
            \bignorm{(T_{-\sum_{n=1}^k 2^{-n}})g}_{M^{p_k}(\R)} &\approx \biggnorm{\Bigparen{\bigip{T_{-\sum_{n=1}^k 2^{-n}})g}{M_ mT_{-2^{-(k+1)}j}\chi_{[0,2^{-(k+1)}]}}}_{m,j\in \Z}}_{\ell^{p_k}(\Z^2)}\\
            &\Ge \biggnorm{\Bigparen{\bigip{T_{-\sum_{n=1}^k 2^{-n}})g}{M_ m\chi_{[0,2^{-(k+1)}]}}}_{m\in\Z}}_{\ell^{p_k}(\Z)}\\
            &\approx \norm{h\cdot\chi_{[\sum_{n=1}^k 2^{-n},\sum_{n=1}^{k+1} 2^{-n}]}}_{M^{p_k}(\R)}=\infty\\  
            &\Eq \infty.
        \end{split}
    \end{align*}
Thus, $g\notin M^{p_k}(\R)$ for any $k\in \N.$   
\end{proof}
\begin{remark}
    It is known that if $\Gc(g,\alpha,\beta)$ is a frame for $L^2(\R)$ and $\alpha \beta =1$, then $\Gc(g,\alpha,\beta)$ is a Riesz basis for $L^2(\R).$ Consequently, for any $1<q\leq 2$ it follows from Theorem \ref{counter_example} that there exist Gabor Riesz bases for which the atoms are in $M^q(\R)$ but not in $M^p(\R)$ for any $1\leq p<q$. Here Riesz bases are Schauder bases for $L^2(\R)$ that are topologically isomorphic to orthonormal bases. However, it is not known whether there exists a Gabor orthonormal basis for which the atom is not in $M^p(\R)$ for any $1\leq p<2.$ Also, we do not know whether there exists a Parserval frame for $L^2(\R)$ whose atom is not in $M^p(\R)$ for any $1\leq p<2$ when $\frac{1}{2}< \alpha\beta\leq 1. $
\end{remark}

\section{Convergence of Gabor Expansions in Modulation Spaces}
In this section, we consider two questions motivated by Theorem \ref{conv_Gabor_modu_original}, Theorem \ref{conv_Gabor_ext} and \cite{HY23}. We introduce some terminology and background first.

Let $X$ be a Banach space with dual space $X^*$. A \emph{Schauder frame} for $X$ is a pair of sequences $\bigparen{(x_n)_{\inN}, (x_n^*)_{\inN}} \subseteq X\times X^*$ for which 
\begin{equation}
\label{Schau_eq}
   x=\sumli \ip{x}{x_n^*}\,x_n, 
\end{equation}
where the series converges in the norm of $X.$ In general, the series in Equation (\ref{Schau_eq}) could be order-sensitive for some $x\in X$. That is, the series associated with some $x\in X$ might become divergent if we permute the order the summation. When the series in Equation (\ref{Schau_eq}) converge unconditionally for all $x\in X$, then we say that $\bigparen{(x_n)_{\inN}, (x_n^*)_{\inN}} \subseteq X\times X^*$ is an \emph{unconditional Schauder frame}. Whenever $\set{x_n^*}_{\inN} \subseteq X^*$ is a sequence such that Equation (\ref{Schau_eq}) holds for all $x\in X$, then we call $\set{x_n^*}_{\inN}$ an \emph{alternative dual} of $\set{x_n}_{\inN}$ and Equation (\ref{Schau_eq}) is called the \emph{reconstruction formula} of $x$. If $\ip{x_n}{x_m^*}=\delta_{mn}$ for all $m,n\in \N$ and the reconstruction formula is unique for every $x\in X$, then we call $\set{x_n}_{\inN}$ a \emph{Schauder basis} for $X$. A sequence $\set{x_n}_{\inN}$ is said to be an \emph{unconditional basis} for $X$ if it is a Schauder basis for $X$ and the associated reconstruction formula converges unconditionally for all $x\in X$.

Let $\Gc(g,\gamma,\alpha,\beta)$ denote the pair of Gabor systems $\bigparen{(M_{\beta n}T_{\alpha k}g)_{k,n\in \Z},(M_{\beta n}T_{\alpha k}g)_{k,n\in \Z}}$. Now, using the language of Schauder frames, Theorem \ref{conv_Gabor_modu_original}(b) means that we can not only extend the reconstruction formula associated with a Gabor Schauder frame $\Gc(g,\gamma,\alpha,\beta)$ from $L^2(\R)$ to $M^{p}(\R)$ but also endow the reconstruction formula with unconditional convergence in $M^{p}(\R)$ for all $1\leq p<\infty$ if the window functions $g$ and $\gamma$ are chosen from $M^1(\R)$. A follow-up question motivated by the extension of Theorem \ref{conv_Gabor_modu_original} is formulated as follows:

\begin{question}
\label{question_51}
Let $g,\gamma\in M^{p}(\R)$ for some $1< p \leq 2$.  Assume that $\Gc(g,\gamma,\alpha,\beta)$ is a Schauder frame for $L^2(\R)$ (with respect to some ordering). Do we still have that \begin{equation}
\label{Gabor_expan_followup}
f = \dousum \ip{f}{\tfshift \gamma} \tfshift g 
\end{equation}
with unconditional convergence of the series for all $f\in M^{q}(\R)$ and $1\leq p\leq q\leq p'<\infty$? 
\end{question}
Note that since $(M^p(\R))^*=M^{p'}(\R)$, we should not expect that Equation (\ref{Gabor_expan_followup}) can be extended to $M^q(\R)$ for some $q>p'.$ Due to Theorem \ref{conv_Gabor_ext}, we see that the answer to Question (\ref{question_51}) is true for some family of functions for all $1<p\leq 2$, for example, $\Gc(\chi_{[0,1]},\chi_{[0,1]},1,1)$. However, it was shown by Heil and Powell in \cite{HP06} that there exist $g,\gamma\in M^{2}(\R)$ for which $\Gc(g,\gamma,1,1)$ is a Schauder basis for $L^2(\R)$ but the associated reconstruction formula does not converge unconditionally for all $f\in L^2(\R)$. We mention that it is still open whether there exists a counterexample to Question (\ref{question_51}) for $1<p<2.$ Specifically, for each $1<p<2$ do there exist $g,\gamma\in M^{p}(\R)$ for which $\Gc(g,\gamma,1,1)$ is a Schauder frame for $L^2(\R)$ but the associated reconstruction formula does not converge unconditionally for all $f\in M^{q}(\R)$ for some $1\leq p\leq q\leq p'<\infty?$ While the answer to Question (\ref{question_51}) is still unknown, we prove several equivalent statements for the answer to be true in the following. 

Recall that Khintchine's Inequality states that for any $1\leq p<\infty$ there exists positive constants $A_p,B_p$ such that for every $N\in \N$ and scalars $c_1,\dots,c_N$ 
$$A_p\norm{(c_n)_{n=1}^N}_{\ell^2} \Le \bignorm{\sum_{n=1}^N c_nR_n}_{L^p([0,1])}\Le B_p\norm{(c_n)_{n=1}^N}_{\ell^2},$$
where $\set{R_n(x)}_{n\geq 0}$ denotes the Rademacher system, which for $n\geq 0$ is defined by $R_n(x)= \text{sign}(\sin(2^n \pi x)).$  Also, recall that a series $\sumli x_n$ is \emph{weakly unconditionally convergent} in $X$ if there exists some constant $C>0$ such that $\sumli |\ip{x_n}{x^*}|\leq C\norm{x^*}_{X^*}$ for all $x^*\in X^*.$

\begin{proposition}
\label{equiva_state_question2}
     Assume that $g,\gamma \in M^{p}(\R)\setminus M^{1}(\R)$ for some $1< p\leq 2$ are such that $\bigparen{\set{M_{\beta n}T_{\alpha k}g}_{k,n\in \Z},\set{M_{\beta n}T_{\alpha k}\gamma}_{k,n\in \Z}}$ is a Schauder frame for $L^2(\R)$ (with respect to some ordering). Then the following statements are equivalent.
     \begin{enumerate}
     \setlength\itemsep{0.8em}
         \item [\textup{(a)}] $\bigparen{\set{M_{\beta n}T_{\alpha k}g}_{k,n\in \Z},\set{M_{\beta n}T_{\alpha k}\gamma}_{k,n\in \Z}}$ is an unconditional Schauder frame for $M^q(\R)$ for any $p\leq q\leq p'$.

             \item [\textup{(b)}] $\bigparen{\set{M_{\beta n}T_{\alpha k}\gamma}_{k,n\in \Z},\set{M_{\beta n}T_{\alpha k}g}_{k,n\in \Z}}$ is an unconditional Schauder frame for $M^q(\R)$ for any $p\leq q\leq p'$.    

  \item [\textup{(c)}] The series $\sum_{n,k\in \Z} \,\ip{f}{M_{\beta n}T_{\alpha k}\gamma}  \,M_{\beta n}T_{\alpha k}g$ is weakly unconditionally convergent in $M^q(\R)$ for all $f\in M^{q}(\R)$ for any $p\leq q\leq p'$.

  \item [\textup{(d)}] The series $\sum_{n,k\in \Z} \,\ip{f}{M_{\beta n}T_{\alpha k}g}\,  M_{\beta n}T_{\alpha k}\gamma$ is weakly unconditionally convergent in $M^q(\R)$ for all $f\in M^{q}(\R)$ for any $p\leq q\leq p'$. 
             
         \item [\textup{(e)}] The analysis operators associated with $g$ and $\gamma$ are bounded linear operators from $M^q(\R)$ to $\ell^q(\Z^2)$ for any $p\leq q\leq p'$. Consequently, there exist some positive constants $A,B$ such that the norm equivalences $$A\,\norm{f}_{M^q(\R)} \leq \bignorm{\bigparen{\ip{f}{M_{\beta n}T_{\alpha k}g}}_{n,k\in \Z}}_{\ell^q(\Z^2)}\leq B\,\norm{f}_{M^{q}(\R)}$$
         and
$$A\,\norm{f}_{M^q} \leq \bignorm{\bigparen{\ip{f}{M_{\beta n}T_{\alpha k}\gamma}}_{n,k\in \Z}}_{\ell^q(\Z^2)}\leq B\,\norm{f}_{M^{q}(\R)}$$
 hold for $M^q(\R)$ for any $p\leq q\leq p'$.
     \end{enumerate}
\end{proposition}   
\begin{proof}
By Lemma \ref{tf_lemmas} (c), $M^q(\R)$ is reflexive for any $p\leq q\leq p'.$ Moreover, by considering $g=\chi_{[0,1]}$ and $\alpha=1=\beta$ in Theorem \ref{conv_Gabor_ext}, we see that $M^p(\R)$ admits an unconditional basis (see also \cite[Theorem 12.3.7]{Gro01}). Consequently, $M^q(\R)$ does not contain a topologically isomorphic copy of $c_0(\Z)$ or $\ell^1(\Z)$ (see \cite[Remark 3.2]{CLS11}).

    (a) $\Leftrightarrow$ (b) It was shown in \cite[Corollary 3.3]{CLS11} that if $\Gc(g,\gamma,\alpha,\beta)$ is an unconditional Schauder frame for a Banach space $X$, then $\Gc(\gamma,g,\alpha,\beta)$ is an unconditional Schauder frame for $X^*$ if $X$ does not contain an topologically isomorphic copy of $\ell^1(\Z).$ This direction now simply follows by duality.
\smallskip

(b) $\Leftrightarrow$ (c) The direction from (b) to (c) is clear. The implication from (c) to (b) is due to the fact that the notions of weakly unconditionally convergent series and unconditionally convergent series are equivalent in Banach spaces that does not contain an isomorphic copy of $c_0$ (for example, see \cite[II.D.Proposition 5]{WO91}). Moreover, the same argument shows the equivalence of (a) and (d).  
\smallskip

(e) $\Rightarrow$ (a) By duality, we see that the synthesis operator $R_\gamma$ associated with $\Gc(\gamma,\alpha,\beta)$ is a bounded linear operator from $M^q(\R)$ to $\ell^q(\Z^2)$ for all $p\leq q\leq p'.$ The unconditional convergence of the reconstruction formula associated with $\Gc(g,\gamma,\alpha,\beta)$ then follows from a similar argument to Lemma \ref{bdd_syn}. Using Inequality (\ref{operator_ineq}), we obtain the norm equivalences. 
\smallskip 

(a) $\Rightarrow$ (e) The proof of this implication closely follows \cite{PY24}. Fix $\psi \in S(\R)$ with $\norm{\psi}_{L^2(\R)}=1$ and let $S=\bigset{V_\psi f\,|\,f\in M^{q}(\R)}$ be equipped with the norm $\norm{\cdot}_{L^q(\R^2)}$. It is not hard to see that $S$ is a Banach space under this norm. By Lemma \ref{tf_lemmas} (e), we see that for any $V_\psi f \in S$ 
    \begin{align*}
        \begin{split}  
    V_\psi f &\Eq V_\psi \Bigparen{\sum_{k,n\in \Z} \ip{f}{M_{\beta n}T_{\alpha k}\gamma} M_{\beta n}T_{\alpha k}g}\\ &\Eq \sum_{k,n\in \Z} \ip{f}{M_{\beta n}T_{\alpha k}\gamma} V_\psi (M_{\beta n}T_{\alpha k}g)\\
    &\Eq 
     \sum_{k,n\in \Z} \ip{V_\psi f}{V_\psi(M_{\beta n}T_{\alpha k}\gamma)} V_\psi (M_{\beta n}T_{\alpha k}g)
     \end{split}
    \end{align*}
    with unconditional convergence of the series. 

Since $\Gc(g,\gamma,\alpha,\beta)$ is an unconditional Schauder frame for $M^q(\R)$, we can enumerate $\Gc(g,\gamma,\alpha,\beta)$ as $\bigset{(g_j)_{j\in\N},  (\gamma_j)_{j\in \N}}$. We then compute 
\begin{equation*}
    \begin{split}
        \bignorm{\bignorm{\bigparen{\ip{V_\psi f}{\gamma_j}\,g_j}_{j\in \N}}_{\ell^2(\N)}}_{L^q(\R^2)} &\Eq \lim_{N\rightarrow \infty} \Bignorm{\Bigparen{\sum_{j=1}^N \,\bigabs{\ip{V_\psi f}{\gamma_j}\,g_j(\cdot)}^2}^{1/2}}_{L^q(\R^2)}\\
        &\Le \lim_{N\rightarrow \infty} k_q^{-1}\Bignorm{\Bignorm{\sum_{j=1}^N\ip{V_\psi f}{\gamma_j}\,g_jR_j(\cdot)}_{L^q([0,1])}}_{L^q(\R^2)}\\
        &\Eq \lim_{N\rightarrow \infty}k_p^{-1}\Bignorm{\Bignorm{\sum_{j=1}^N\ip{V_\psi f}{\gamma_j}\,g_j(\cdot )R_j}_{L^q(\R^2)}}_{L^q([0,1])}\\
    \end{split}
\end{equation*}
where the second inequality is a direct application of Khinchine's Inequalities. 

By \cite[Lemma 2.1]{LT23}, there exists a constant $K>0$ such that for any $|\epsilon_{k,n}|\leq 1,$ 
$$\Bignorm{\displaystyle \sum_{k,n\in \Z} \ip{V_\psi f}{V_\psi(M_{\beta n}T_{\alpha k}\gamma)} \epsilon_{k,n} V_\psi (M_{\beta n}T_{\alpha k}g)}_{L^q(\R^2)}\Le K\norm{V_\psi f}_{L^q(\R^2)}.$$

Consequently, 
\begin{equation*}
    \begin{split}
       \lim_{N\rightarrow \infty}k_p^{-1}\Bignorm{\Bignorm{\sum_{j=1}^N\ip{V_\psi f}{\gamma_j}\,g_j(\cdot )R_j}_{L^q(\R^2)}}_{L^q([0,1])} & \Le Kk_p^{-1} \norm{\norm{V_\psi f}_{L^q(\R^2)}}_{L^q([0,1])}\\
        &\Eq Kk_p^{-1} \norm{f}_{M^q(\R)}.
    \end{split}
\end{equation*}

On the other hand, without loss of generality, we assume that $\norm{V_\psi \gamma}_{L^q([0,\alpha]\times [0,\beta])}>0.$ Let $I_{\ell,m}$ denote $[\alpha \ell, \alpha(\ell+1)]\times [\beta m, \beta(m+1)]$. Then a straightforward computation yields
\begin{align*}
    \begin{split}
    \bignorm{\bignorm{\bigparen{\ip{V_\psi f}{g_j}\,\gamma_j(s)}}_{\ell^2(\N)}}^q_{L^q(\R^2)}&\Eq \sum_{i,j\in \Z}\int_{I_{i,j}} \bignorm{\bigparen{\ip{f}{M_{\beta n}T_{\alpha k}g} V_\psi(M_{\beta n}T_{\alpha k}\gamma)}_{k,n\in \Z}}_{\ell^2(\Z^2)}^q \\
    &\Ge \sum_{\ell,m\in \Z}\int_{I_{\ell,m}}\bigabs{\ip{f}{M_{\beta m}T_{\alpha \ell}g}V_\psi \gamma(x-\alpha \ell,w-\beta m)}^{q}\\
    &\Eq \sum_{\ell,m\in \Z} \bigabs{\ip{f}{M_{\beta m}T_{\alpha \ell}g}}^q\int_{I_{0,0}} \bigabs{V_\psi \gamma(x, w)}^{q} \,dxdw\\
    &\Eq C\,\bignorm{\bigparen{\ip{f}{M_{\beta m}T_{\alpha \ell}g}}_{\ell,m\in \Z}}^q_{\ell^q(\Z^2)}.
     \end{split}
\end{align*}
Thus, $C_g\colon M^{q}(\R)\rightarrow \ell^q(\Z^2)$ is a bounded linear operator.
Using the equivalence of (a) and (b) and a similar argument used to show the boundedness of $C_g$, we obtain the boundedness of $C_\gamma.$ The boundedness of the synthesis operators associated with $g$ and $\gamma$ then follow by duality. The norm equivalences now follow from Inequality (\ref{operator_ineq}).\qedhere
\end{proof}

\begin{remark}
We remark that if $g,\gamma \in M^{p,q}(\R)\setminus M^{1,1}(\R)$ for some $1< p,q\leq 2$, then the following statements still hold by using the same proof of Proposition \ref{equiva_state_question2}. 
 \begin{enumerate}
     \setlength\itemsep{0.8em}
         \item [\textup{(a)}] $\bigparen{\set{M_{\beta n}T_{\alpha k}g}_{k,n\in \Z},\set{M_{\beta n}T_{\alpha k}\gamma}_{k,n\in \Z}}$ is an unconditional Schauder frame for $M^{p_1,q_1}(\R)$ for any $p\leq p_1\leq p'$ and $q\leq q_1\leq q'$.

             \item [\textup{(b)}] $\bigparen{\set{M_{\beta n}T_{\alpha k}\gamma}_{k,n\in \Z},\set{M_{\beta n}T_{\alpha k}g}_{k,n\in \Z}}$ is an unconditional Schauder frame for $M^{p_1,q_1}(\R)$ for any $p\leq p_1\leq p'$ and $q\leq q_1\leq q'$.

  \item [\textup{(c)}] The series $\sum_{n,k\in \Z} \,\ip{f}{M_{\beta n}T_{\alpha k}\gamma}  \,M_{\beta n}T_{\alpha k}g$ is weakly unconditionally convergent in $M^{p_1,q_1}(\R)$ for all $f\in M^{p_1,q_1}(\R)$ and any $p\leq p_1\leq p'$ and $q\leq q_1\leq q'$.

  \item [\textup{(d)}] The series $\sum_{n,k\in \Z} \,\ip{f}{M_{\beta n}T_{\alpha k}g}\,  M_{\beta n}T_{\alpha k}\gamma$ is weakly unconditionally convergent in $M^{p_1,q_1}(\R)$ for all $f\in M^{p_1,q_1}(\R)$ and any $p\leq p_1\leq p'$ and $q\leq q_1\leq q'$.
  \end{enumerate}
\end{remark}

We close the discussion on Question (\ref{question_51}) by showing that although it is not clear whether $ f=\sum_{k,n\in \Z}\ip{f}{\tfshift \gamma} \tfshift g $ holds with the unconditional convergence of the series in the same space where $g,\gamma$ are, we can have at least that $f=\dousum\, \ip{f}{\tfshift \gamma} \tfshift g $ holds with the unconditional convergence of the series in some modulation space that is bigger than the space where $g$ and $\gamma$ are.

\begin{definition}	
	 A pair of real numbers $(p,p_1)$ is an \emph{extensible pair} if $(p,p_1) \in [1,2] \times [1,\infty)$ and $p_1< \dfrac{p}{2p-2}.$
\end{definition}

A straightforward computation shows that if $(p,p_1)$ is an extensible pair of real numbers if and only if $1\leq \frac{pp_1}{p+p_1-pp_1}, \frac{pp_1}{p+2p_1-2pp_1}<\infty.$
We will prove the following result. 
\begin{theorem}
	\label{exten_conv_in_modu}
	Fix $\alpha,\beta >0$ and let $(p,p_1)$ be an extensible pair. Assume that $g,\gamma \in \mpq$ are such that $\bigparen{\set{M_{\beta n}T_{\alpha k}g}_{k,n\in \Z},\set{M_{\beta n}T_{\alpha k}\gamma}_{n,k\in \Z}}$ is a Schauder frame for $L^2(\R)$. Then 
	$$f\Eq\sum_{k,n\in \Z} \ip{f}{M_{\beta n}T_{\alpha k}\gamma}  M_{\beta n}T_{\alpha k}g $$
	with the unconditional convergence of the series in $M^{\frac{pp_1}{p+2p_1-2pp_1}}(\R)$ for all $f\in \mpqone$.
	\end{theorem}

To prove Theorem \ref{exten_conv_in_modu}, it suffices to show that the $C_g\colon M^{p_1}(\R)\rightarrow \ell^{\frac{pp_1}{p+p_1-pp_1}}(\Z^2)$ and $R_\gamma\colon \ell^{\frac{pp_1}{p+p_1-pp_1}}(\Z^2)\rightarrow M^{\frac{pp_1}{p+2p_1-2pp_1}}(\R)$ are bounded linear operators. Then the conclusion follows from the same argument used in Theorem \ref{conv_Gabor_ext}.
We will need several lemmas below to prove the boundedness of these two operators.
\begin{lemma}
	\label{STFT_conti}
	Let $(p,p_1)$ be an extensible pair. Assume that $g\in \mpq$ and $f\in \mpqone$. Then the following statements hold\begin{enumerate}
 \setlength\itemsep{0.8em}
	\item [\textup{(a)}]$V_{g} f(x,w)$ is continuous in $\mathbb{R}^{2}$. 
	\item [\textup{(b)}] $|V_{g} f|\Le \parenspace{|V_{\gamma} f|\ast |V_g \gamma|}$ for any $\gamma\in S(\R)$ with $\norm{\gamma}_{L^2(\R)}=1$.	\end{enumerate}
\end{lemma}
\begin{proof}
	(a)
	Let $\set{g_n}_{\inN}\subseteq S(\R)$ be a sequence that converges to $g$ in $\mpq$ and let $\set{f_n}_{\inN} \subseteq S(\R)$ be a sequence that converges to $f $ in $\mpqone$. Fix $\gamma\in S(\R)$ with $\norm{\gamma}_{L^2(\R)}=1$. Then by \cite[Lemma 11.3.3]{Gro01}, we have 
	\begin{equation}
	\label{convo_inequa_extend}
	|V_{g_n}f_n |\Le \bigparen{|V_{\gamma} f_n|\ast |V_{g_n}\gamma| } \Le \bigparen{|V_{\gamma}f_n|\ast |V_{g_n}\gamma| \ast |V_{\gamma}\gamma|}.
	\end{equation}
	By Young's convolution inequality, we see that 
 \begin{equation}
 \label{exten_eq1}
    \begin{split}
 \norm{V_{g_n}f_n }_{L^\infty(\R^2)} &\Le \bignorm{V_{\gamma}f_n}_{L^{p_1}(\R^2)}\,\,\bignorm{\,|V_{g_n}\gamma| \ast |V_{\gamma}\gamma|\,}_{L^{p_1'}(\R^2)}\\&\Eq \norm{f_n}_{\mpqone}\,\,\bignorm{\,|V_{g_n}\gamma| \ast |V_{\gamma}\gamma|\,}_{L^{p_1'}(\R^2)},
   \end{split}
 \end{equation}
 where $p_1'$ is the conjugate exponent of $p_1$.
Since $1+\frac{1}{p'}-\frac{1}{p}= 2-\frac{2}{p} \leq 1$, we have that $\frac{1}{2-\frac{2}{p}}\Eq \frac{p}{2p-2}\Ge 1$. 
Applying Young's convolution inequality again, we obtain 
\begin{equation}
\label{exten_eq2}
\bignorm{\,|V_{g_n}\gamma| \ast |V_{\gamma}\gamma|\,}_{L^{p_1'}(\R^2)}\Le \bignorm{V_{g_n}\gamma}_{L^{p}(\R^2)}\norm{V_{\gamma}\gamma}_{L^{\frac{p}{2p-2}}(\R^2)}.
\end{equation}
Combining (\ref{exten_eq1}) and (\ref{exten_eq2}), we reach the estimate
\begin{equation}
\label{L_infinity_norm_STFT}
\norm{V_{g_n}f_n }_\infty \Le  
\norm{g_n}_{\mpq} \,\norm{f_n}_{\mpqone} \,\norm{V_{\gamma}\gamma}_{L^{\frac{p}{2p-2}}(\R^2)}.
\end{equation}
Therefore, $V_{g_n}f_n (x,w)$ is a Cauchy sequence in $L^{\infty}(\R^{2})$. We then define $V_gf(x,w)$ pointwise by $$V_gf(x,w)\Eq \lim_{n\rightarrow \infty} V_{g_n}f_n (x,w).$$
By \cite[Theorem 11.2.3]{Gro01}, $V_{g_n}f_n(x,w)$ is continuous in $\R^2$ for all $n$, and hence $V_gf(x,w)$ is continuous in $\R^2.$ 
So, it remains to show that $V_gf$ is well-defined. Let $\set{h_m}_{m\,\in\,\N}, \set{w_m}_{m\,\in\,\N}\subseteq S(\R)$ be two sequences that converge $f,g$ in $\mpqone$ and $\mpq$, respectively. By the Triangle Inequality, we have $$|V_{g_n}f_n-V_{w_m}h_m| \Le |V_{g_n}(f_n-h_m)|+|V_{g_n-w_m}h_m| .$$
Arguing similarly to the way we obtained Inequality (\ref{L_infinity_norm_STFT}), we see that 
$$\norm{V_{g_n}(f_n-h_m)}_{L^\infty(\R^2)} \Le \norm{g_n}_{\mpq} \,\norm{f_n-h_m}_{\mpqone} \,\norm{V_{\gamma}\gamma}_{L^{\frac{p}{2p-2}}(\R^2)} $$
and 
$$\norm{V_{g_n-w_m}h_m}_{L^\infty(\R^2)} \Le \norm{g_n-w_m}_{\mpq}\, \norm{h_m}_{\mpqone} \,\norm{V_{\gamma}\gamma}_{L^{\frac{p}{2p-2}}(\R^2)}.$$ 
Therefore, $V_gf $ is well-defined.
\medskip

(b) This statement follows from the lower inequality of Inequality (\ref{convo_inequa_extend}) and the fact that $|V_{\gamma} f_n|$ and $|V_{g_n}\gamma|$ converges to $|V_{\gamma} f|$, $|V_{g}\gamma|$  in $L^\infty$-norm, respectively.
\end{proof}



\begin{lemma}
	\label{extended_inclusion_property}
	    Let $(p,p_1)$ be an extensible pair. Assume that $f\in \mpqone$ and $g\in \mpq$. Then $\bigparen{\norm{V_{g} f\cdot \chi_{[k,k+1]\times[n,n+1]}}_{L^\infty(\R)}}_{k,n\in \Z} \in \ell^{\frac{pp_1}{p+p_1-pp_1}}(\Z^2)$ 
    and $$\bignorm{ \bigparen{\norm{V_{g} f\cdot \chi_{[k,k+1]\times[n,n+1]}}_{L^\infty(\R)}}_{k,n\in \Z}   }_{\ell^\frac{pp_1}{p+p_1-pp_1}(\Z^2)}\Le C \norm{g}_{\mpq}\,\norm{f}_{\mpqone}.$$
	    for some constant $C>0$.
	\end{lemma}
\begin{proof}
	Let $\gamma \in S(\R)$ with $\norm{\gamma}_2=1$. By Lemma \ref{STFT_conti} (b), we have $$|V_{g}f|\leq  \parenspace{|V_\gamma f|\ast |V_f\gamma|\ast |V_\gamma\gamma|} .$$
	Using \cite[Theorem 11.1.5]{Gro01} and \cite[Theorem 12.2.1]{Gro01}, we obtain $$\bignorm{ \bigparen{\norm{V_{g} f\cdot \chi_{[k,k+1]\times[n,n+1]}}_{L^\infty(\R)}}_{k,n\in \Z}   }_{\ell^\frac{pp_1}{p+p_1-pp_1}(\Z^2)}\Le C\bignorm{|V_\phi f|\ast |V_g\phi|}_{L^{\frac{pp_1}{p+p_1-pp_1}}(\R^{2})},$$
	for some $C>0$. The statement then follows from Young's convolution inequality.
	\end{proof}

We are now ready to prove the boundedness of the analysis operator and the synthesis operator associated with a Gabor system $\Gc(g,\alpha,\beta)$ when the atom $g$ belongs to $M^p(\R)$ for some  $1\leq p\leq 2$.  We mention that the proofs of following two lemmas follow closely to \cite[Theorem 12.2.3]{Gro01} and \cite[Theorem 12.2.4]{Gro01}.
\begin{lemma}
\label{bdd_ana_large_modu}
Fix $\alpha, \beta>0$.
	and let $(p,p_1)$ be an extensible pair. Then for any $g\in \mpq$ the analysis operator $C_{g}$ associated with $\Gc(g,\alpha,\beta)$
	is a bounded linear operator from $M^{p_1}(\R)$ to $\ell^{\frac{pp_1}{p+p_1-pp_1}}(\Z^{2})$. 
	\end{lemma}
\begin{proof}
	Note that $\ip{f}{M_{\beta n}T_{\alpha_k}g}= V_{g}f (\alpha k,\beta n)$ is well-defined by Lemma \ref{STFT_conti} (a). By \cite[Prop. 11.1.4]{Gro01} and Lemma \ref{extended_inclusion_property}, we obtain
	\begin{align*}
	\begin{split}
	\bignorm{C_g(f)}_{\ell^{\frac{pp_1}{p+p_1-pp_1}, \frac{qq_1}{q+q_1-qq_1}}(\Z^{2})} &\Eq \bignorm{\bigparen{V_{g}f (\alpha k,\beta n)}_{k,n\in \Z}}_{\ell^{\frac{pp_1}{p+p_1-pp_1}}(\Z^{2})}\\ &\Le \bignorm{ \bigparen{\norm{V_{g} f\cdot \chi_{[k,k+1]\times[n,n+1]}}_{L^\infty(\R)}}_{k,n\in \Z}   }_{\ell^\frac{pp_1}{p+p_1-pp_1}(\Z^2)} \\& \Le C\, \norm{g}_{\mpq}\,\norm{f}_{\mpqone}. \qedhere
	\end{split}
	\end{align*}
	\end{proof}
\begin{lemma}
\label{bdd_syn_large_modu}
		Fix $\alpha, \beta>0$ and let $(p,p_1)$ be an extensible pair. Then for any $g\in \mpq$ the synthesis operator $R_{g}$ associated with $\Gc(g,\alpha,\beta)$ is a bounded linear operator from $\ell^{p_1}(\Z^{2})$ to $M^{\frac{pp_1}{p+p_1-pp_1}}(\R)$. \\
		Consequently, the series $\sum_{k,n\,\in\, \mathbb{Z}^2} c_{kn}M_{\beta n}T_{\alpha k}g $ converges unconditionally in $M^{\frac{pp_1}{p+p_1-pp_1}}(\R)$ for any $(c_{kn})_{k,n\in \Z}\in \ell^{p_1}(\Z^2).$
	\end{lemma}
\begin{proof}
We will first show that $R_g(c)\in S'(\R)$ for any $c=(c_{kn})_{k,n\in \Z}\in \ell^{p_1}(\Z^2)$.
Note that since $(p,p_1)$ is an extensible pair, $M^p(\R)\subseteq M^{p_1'}(\R).$ Then an argument similar to the one used in the beginning of the proof of Theorem \ref{bdd_syn} shows that $R_g(c)\in S'(\R)$.
	Next, fix $\gamma \in S(\R)$ and let $G= |V_\gamma g|$. Then we have
 $$|V_\gamma \bigparen{R_g (c)}|\Le \sum_{k,n\,\in\, \mathbb{Z}^d} |c_{kn}| \,|V_{\gamma} (M_{\beta_n} T_{\alpha_k} g)|\Le  \sum_{k,n\,\in\, \mathbb{Z}^d} |c_{kn}|\, G(x-\alpha k, w-\beta n).$$
	For each $m,\ell \in \Z$ we define $$a_{m\ell} \Eq \sup_{(x,w)\,\in\, Q_\alpha \times Q_\beta} G(x-\alpha m, w-\beta \ell),$$
 where $Q_\alpha = [0,\alpha].$
	Since $G(x,w)\leq \sum_{m,\ell\in \Z^d} a_{m\ell} \,T_{(\alpha m, \beta \ell)} \,\chi_{Q_{\alpha}\times Q_{\beta}},$ we compute 
	\begin{align*}
	 \bigabs{V_\gamma \bigparen{R_g (c)}}&\Le  \sum_{k,n,m,\ell\,\in\, \mathbb{Z}} |c_{kn}| |a_{m,\ell}| \,T_{(\alpha (k+m), \beta (n+\ell))}\,\chi_{Q_{\alpha}\times Q_{\beta}} \quad \text{( Let $i=k+m$ and $j=n+\ell.$ )}\\
	 &\Eq  \sum_{i,j\,\in\, \mathbb{Z}}\Bigparen{\sum_{k,n\,\in\, \mathbb{Z}} |c_{kn}| \,|a_{i-k,j-n}| }\,T_{(\alpha i, \beta j)}\,\chi_{Q_{\alpha}\times Q_{\beta}}\\
	 & \Eq \sum_{i,j\,\in\, \mathbb{Z}}\bigparen{|c|\ast |a|}\,T_{(\alpha i, \beta j)}\,\chi_{Q_{\alpha}\times Q_{\beta}}.
	 \end{align*} 
 Consequently, 
 \begin{align*}
	\bigabs{V_\gamma \bigparen{R_g (c)}}_{L^{\frac{pp_1}{p+p_1-pp_1}}(\R^{2})} &\Le C\, \bignorm{|c|\ast |a|}_{\ell^{\frac{pp_1}{p+p_1-pp_1}}(\Z^{2})}.
	\end{align*}
 Since $(a_{m,\ell})_{m,\ell\,\in\, \Z}\in \ell^{p}(\Z^{2})$ by Lemma \ref{extended_inclusion_property}, we reach the following estimate
	\begin{align*}
	\norm{V_\gamma \bigparen{R_g (c)}}_{L^{\frac{pp_1}{p+p_1-pp_1}}(\R^{2})} &\Le C\, \bignorm{|c|\ast |a|}_{\ell^{\frac{pp_1}{p+p_1-pp_1}}(\Z^{2})}\\
	&\Le C \norm{a}_{\ell^{p}(\Z^{2})} \norm{c}_{\ell^{p_1}(\Z^2)}
	\end{align*}
by Young's convolution inequality. Finally, arguing similarly to Lemma \ref{bdd_syn}, we obtain the unconditional convergence of $R_g(c)$ for any $c\in \ell^{p_1}(\Z^2).$  \qedhere
	\end{proof}
\begin{remark}

	(a) It is not necessary to select $g$ and $\gamma$ in Theorem \ref{exten_conv_in_modu} from the same modulation space. For example, if we select $\gamma$ from $M^{q}(\R)$, then an additional assumption that $(q,\frac{pp_1}{p+p_1-pp_1})$ is an extensible pair is required. Then we will obtain that $$f\Eq\sum_{n,k\in \Z} \ip{f}{M_{\beta n}T_{\alpha k}\gamma}  M_{bn}T_{ak}g $$ with the unconditional convergence in $M^{\frac{pp_1+q_1+qp_1-2qpp_1}{pqp_1}}(\R)$ for all $f\in \mpqone.$
	\smallskip
	
	(b) We obtain Thereom \ref{conv_Gabor_modu_original} by letting $p=1$ in Theorem \ref{exten_conv_in_modu}. However, if $1<p\leq 2$, then $p\leq \frac{pp_1}{p+2p_1-2pp_1}$. Therefore, we only obtain convergence of the series in the Equation (\ref{Gabor_expan_followup}) in some larger modulation space. \qeddef
	\end{remark}


The next question is about the \emph{classification of alternative duals} that was studied in \cite{HY23} for frames in Hilbert spaces. 
Let $\Gc(g,\gamma,\alpha,\beta)$ be a Schauder frame for $M^{p,q}(\R)$ for some $1\leq p,q\leq 2.$ It is known that a Gabor system $\set{\tfshift g}_{n,k\in \Z}$ could admit more than one alternative dual. As a result, it is possible that we might lose the unconditional convergence of the reconstruction formula when we shift from one alternative dual to another. If we use the unconditional convergence of the reconstruction formula to classify ``good" and "bad" alternative duals, then
 Theorem \ref{conv_Gabor_modu_original} and Theorem \ref{conv_Gabor_ext} provide some conditions for an alternative dual associated with a Gabor system to be a good alternative dual. Specifically, if $\set{\tfshift g}_{n,k\in \Z}$ is a Gabor system with $g\in C^1_{AC}[a,b]$, then every alternative dual of the form $\set{\tfshift \gamma}_{n,k\in \Z}$ with $\gamma\in C^1_{AC}[a,b]$ is a good alternative dual by Theorem \ref{conv_Gabor_ext}. However, it is not necessary that every alternative dual of $\set{\tfshift g}_{n,k\in \Z}$ is in the form of a Gabor system with the same translation and modulation parameter. Even if an alternative dual is of the form  $\set{\tfshift \gamma}_{n,k\in \Z}$ for some $\gamma$, the atom needs not to be in the same modulation space as $g$. We provide an example below where the window functions of a Gabor Schauder frame are not in the same modulation space.
 \begin{example}
      Let $g(x)=\frac{x}{2}\cdot \chi_{[0,1)}+(1-\frac{x}{2}) \cdot \chi_{[1,2)}.$ By \cite[Proposition 12.1.6]{Gro01}, we see that $g\in M^1(\R).$ Next, let $h=\frac{1}{2} \chi_{[0,2]}.$ We mentioned before that $h$ is not in $M^1(\R)$ by \ref{bdd_of_gchi} and the fact the $M^p(\R)$ is closed under the Fourier transform for any $1\leq p\leq \infty.$ Using Theorem \ref{painless_non_expan}, we obtain that both $\Gc(g,1,\frac{1}{2})$ and  $\Gc(h,1,\frac{1}{2})$ are frames for $L^2(\R)$ (see also \cite[Corollary 11.7.1]{Chr16}). We will show that $\Gc(g,1,\frac{1}{2})$ and $\Gc(h,1,\frac{1}{2})$ are alternative dual of each other.
        By \cite[Theorem 12.3.4]{Chr16}, it suffices to show that $$\displaystyle\sum_{k\in\Z} g(x-k-2n)\,h(x-k)= \frac{1}{2}\delta_{n,0} \qquad \textup{for a.e. } x\in [0,1]. $$
        Both $g$ and $h$ are supported in $[0,2]$, so $\sum_{k\in\Z} g(x-k-2n)h(x-k)=0$ if $n\neq 0.$
       If $n= 0$, then $\sum_{k\in\Z} g(x-k)h(x-k)=\frac{1}{2}$ for a.e. $x\in [0,1]$ because $\sum_{k\in \Z} g(x-k) = 1$  a.e..
 \end{example}
Alternative duals associated with a given Gabor system can be very different from each other. 
Consequently, unconditional convergence of the reconstruction formula associated with a Gabor Schauder frame becomes much more uncertain. We accordingly formulate the following question.
\begin{question} \label{question_52}
Fix $1\leq p<\infty.$ Assume that $\bigparen{\set{M_{\beta n}T_{\alpha k}g}_{k,n\in \Z},\set{M_{\beta n}T_{\alpha k}\gamma}_{k,n\in \Z}} \subseteq M^p(\R) \times (M^p(\R))^*$ is a Schauder frame for $M^p(\R).$
    Then for what Gabor systems $\Gc(g,\alpha,\beta)\subseteq M^p(\R)$ do we have $$f\Eq\sum_{k,n\in \Z} \ip{f}{\gamma_{k,n}}  M_{bn}T_{ak}g $$
    with the unconditional convergence of the series for all $f\in M^p(\R)$ and every alternative dual $\set{\gamma_{n,k}}_{n,k\in \Z}\subseteq M^{p'}(\R)$?
\end{question}
 It was shown by Heil and the author in \cite{HY23} that the only possibility for a Gabor frame in $L^2(\R)$ to admit ``good alternative duals" only is that the Gabor frame must be a \emph{Riesz basis} for $L^2(\R)$ plus at most finitely many elements. Here Riesz bases are Schauder bases that are topologically isomorphic to orthonormal bases for $L^2(\R).$ We close this paper by presenting a complete answer to Question \ref{question_52}. Even more, our result applies to all Banach spaces that does not contain an isomorphic copy of $c_0.$  

We first extend \cite[Theorem 3.4]{HY23} from separable Hilbert spaces to separable Banach spaces. 
\begin{lemma}
\label{const_alter_Bana_spaces}
    Let $\bigparen{(x_n)_{\inN}, (x_n^*)_{\inN}} \subseteq X\times X^*$ be a Schauder frame for $X$. Assume that there exists a sequence of scalars $(c_n)_{\inN}$ for which $\sumli c_nx_n$ converges to some nonzero element $x_0\in X$. Then there exists a sequence $\set{y_n^*}\subseteq X^*$ for which \begin{enumerate}
    \setlength\itemsep{0.8em}
        \item [\textup{(a)}] $\set{y_n^*}_{\inN}$ is an alternative dual of $\set{x_n}\inN,$
        \item [\textup{(b)}] $\ip{x_0}{y_n^*}=c_n~$  for all $n\in \N$.
    \end{enumerate} 
\end{lemma}
\begin{proof}
Let $M=\text{span}\set{x_0}$. Then we define $T\colon
M\rightarrow F$ by $T(x)=T(\alpha x_0)=\alpha.$ By Hahn-Banach Extension Theorem, we extend $T$ to be a bounded linear functional $T\colon X\rightarrow F$. Then the projection operator $P\colon X\rightarrow X$ defined by $P(x)=T(x)x_0$ satisfies that $P^2=P$ and $\text{Range}(P)=M$. Consequently, every element $x$ in $X$ can be uniquely written $x=x_M+x_{\text{ker}(P)}$ for some elements $x_M\in M$ and $x_{\text{ker}(P)}\in \text{ker}(P)$.
Next, for each $n\in\N$ we define $$y_n^* = x_n^*\circ P_{\text{ker}(P)} + \overline{c_n}(T\circ P_M),$$ where $P_{\text{ker}(P)},P_M$ are projections onto $N,M$, respectively. We will show that $\set{y_n^*}_{\inN}$ is the desired alternative dual. Clearly, $\set{y_n^*}_{\inN}\subseteq X^*.$ For any $x\in X$ we have 
 \begin{align*}
 \begin{split}
 \sumli \ip{x}{y_n^*}\,x_n &\Eq \sumli \ip{x}{x_n^*\circ P_{\text{ker}(P)}}\,x_n + \sumli \ip{x}{\overline{c_n}(T\circ P_M)}\,x_n \\
 &\Eq P_{\text{ker}(P)}(x) + T(P_Mx)\sumli c_nx_n  \\
 &\Eq P_{\text{ker}(P)}(x) + T(P_Mx)x_0\\
 &\Eq P_{\text{ker}(P)}(x)+P_M(x)\\
 &\Eq x.
  \end{split}
 \end{align*}  
 Therefore, $\set{y_n^*}_{\inN}$ is an alternative dual of $\set{x_n}_{\inN}$. 
 
 Finally, for each $n\in \N$ we have
 $$\ip{x_0}{y_n^*} \Eq \ip{x_0}{x_n^*\circ P_N} + \ip{x_0}{c_nT\circ P_M} \Eq c_n.$$ \qedhere
\end{proof}

An immediate corollary of Theorem \ref{const_alter_Bana_spaces} is that Schauder bases and Schauder frames can be distinguished by the number of alternative duals they possess.

\begin{corollary} 
\label{uniq_alter_Schauder_basis}
Let $(\set{x_n}_{\inN}, \set{x_n^*}_{\inN})\subseteq X\times X^*$ be a Schauder frame for $X$. Then $\set{x_n}_{\inN}$ is a Schauder basis if and only if it admits a unique alternative dual.
\end{corollary}
\begin{proof}
 If $\set{x_n}_{\inN}$ is a Schauder basis for $X$, then it has an alternative dual that is biorthogonal to $\set{x_n}_{\inN}$. The uniqueness then follows from the biorthogonality. 
 
 Next, suppose to the contrary that $\set{x_n}_{\inN}$ admits a unique alternative dual but is not a Schauder basis for $X.$
 Note that since $\set{x_n}_{\inN}$ admits a unique alternative dual, we have that $x_n\neq 0$ for every $n$. If $\set{x_n}_{\inN}$ is not a Schauder basis, then there must exist some $x_0\in X$ for which there exist two distinct sequences of scalars $(c_n)_{\inN}$ and $(d_n)_{\inN}$ such that $$x_0=\sumli c_nx_n = \sumli d_nx_n.$$
 If $x_0\neq 0$, then by Lemma \ref{const_alter_Bana_spaces}, we obtain two different alternative duals, which is a contradiction.  If $x_0=0$, then for any $c>0$ we have $$cx_1= \sumli c_n x_n+cx_1= \sumli d_n x_n+cx_1.$$

 So, $cx_1\neq 0$ can be expressed in two different ways in terms $x_n$. Applying Lemma \ref{const_alter_Bana_spaces} again, we obtain two different alternative duals of $\set{x_n}_{\inN}$, which is a contradiction.
\end{proof}

We will need the following characterization of Banach spaces that do not contain a topologically isomorphic copy of $c_0$ due to Casazza and Christensen (\cite{CC96}).
\begin{theorem} (\cite[Theorem 3.2]{CC96})
\label{CC_uncond}
Let $X$ be a separable Banach space. The following statements are equivalent.
    \begin{enumerate}
        \item [\textup{(a)}] $X$ does not contain a topologially isomorphic copy of $c_0.$
        \item [\textup{(b)}] Let $\set{x_n}_{\inN} \subseteq X$ be a sequence that does not contain infinitely many zeros. If $\sumli a_nx_n$ converges unconditionally whenever it converges for any sequence of scalars $(a_n)_{\inN}$, then $\set{x_n}_{\inN}$ is an unconditional basis plus at most finitely many elements. 
        \qeddef
    \end{enumerate} 
\end{theorem}

Now we are ready to present the answer to Question (\ref{question_52}).
\begin{theorem}
\label{BS_uncon_alter_nearuc}
    Let $X$ be a separable Banach spaces that does not contain a copy of $c_0$. Assume that $\bigparen{\set{x_n}_{\inN}, \set{x_n^*}_{\inN}}\subseteq X\times X^*$ is a Schauder frame for $X$ and $\set{x_n}_{\inN}$ does not contain infinitely many zeros. Then the following statements are equivalent.
    \begin{enumerate}
    \setlength\itemsep{0.8em}
        \item [\textup{(a)}] $\set{x_n}_{\inN}$ is an unconditional basis for $X$ plus at most finitely many elements.
        \item [\textup{(b)}] $x = \sumli \ip{x}{x_n^*}x_n$ with unconditional convergence of the series in $X$ for every $x\in X$ and every alternative dual $\set{x_n^*}_{\inN}$.
    \end{enumerate}
    \end{theorem}
    \begin{proof}
        ($\Rightarrow$) 
         Let $A=\set{n_j}_{j=1}^N$ be a subset of $\N$ for which $\set{x_n}_{n \notin A}$ is an unconditional basis for $X$ and let $\set{y_n^*}_{n\notin A}$ be the biorthogonal system of $\set{x_n}_{n\notin A}.$ For any $x\in X$, we have 
         \begin{equation}
         \label{expansion_eq1}
          x = \sum_{n\notin A} \ip{x}{x_n^*}\,x_n + \sum_{j=1}^N \ip{x}{x_{n_j}^*}\,x_{n_j}.
         \end{equation}
         Substituting $x_{n_j} = \sum_{n\notin A} \ip{x_{n_j}}{y_n^*}\,x_n$ into Equation (\ref{expansion_eq1}), we obtain 
         $$ x = \sum_{n\notin A} \ip{x}{x_n^*}\,x_n + \sum_{n\notin A} \bigip{x}{\sum_{j=1}^N\ip{x_{n_j}}{y_n^*}x^*_{n_j}}\,x_n.$$
Since $\set{x_n}_{n\notin A}$ is an unconditional basis, we must have $$\ip{x}{y_n^*}= \ip{x}{x_n^*}+\bigip{x}{\sum_{j=1}^N\ip{x_{n_j}}{y_n^*}x^*_{n_j}},$$
which implies \begin{equation}
\label{substi_eq}
    x_n^*= y_n^*- \sum_{j=1}^N\ip{x_{n_j}}{y_n^*}\,x_{n_j}^*,
\end{equation}
for any $n\notin A.$ 

Next, let $\sigma$ be a permutation of $\N.$ Note that since $A$ is a finite subset of $\N$, the series $\sumli \ip{x}{x^*_{\sigma(n)}}x_{\sigma(n)}$ converges if and only if 
          $\sum_{n\in \N, \,\sigma(n)\notin A} \ip{x}{x^*_{\sigma(n)}}x_{\sigma(n)}$ converges. Using Equation (\ref{substi_eq}), we see that \begin{equation*}
          \begin{split}
               \sum_{n\in\N,\, \sigma(n)\notin A}^\infty \ip{x}{x^*_{\sigma(n)}}x_{\sigma(n)} &\Eq \sum_{n\in \N,\, \sigma(n)\notin A} \bigip{x}{y_{\sigma(n)}^*- \sum_{j=1}^N\ip{x_{n_j}}{y_{\sigma(n)}^*}x_{n_j}^*}x_{\sigma(n)}\\
               & \Eq \sum_{n\in \N, \,\sigma(n)\notin A} \biggip{x-\sum_{j=1}^N \ip{x}{x^*_{n_j}}x_{n_j}}{y_{\sigma(n)}^*}\,x_{\sigma(n)},
               \end{split}
               \end{equation*}
        which converges since $\set{x_n}_{n\notin A}$ is an unconditional basis for $X.$
        \medskip
        
($\Leftarrow$) Suppose to the contrary that $\set{x_n}_{\inN}$ is not an unconditional basis plus at most finitely many elements. By Theorem \ref{CC_uncond}, there exists a sequence of scalars $(c_n)_{\inN}$ for which $\sumli c_nx_n$ converges but not unconditionally. Arguing similarly to Corollary \ref{uniq_alter_Schauder_basis}, we can find an alternative dual $\set{y_n^*}_{\inN}$ and an element $x_0\in X$ for which $\ip{x_0}{y_n^*}=c_n$ for all $n$ and $x_0= \sumli \ip{x_0}{y_n^*}x_n$, where the series converges but not unconditionally, which is a contradiction.
    \end{proof}

\begin{corollary}
\label{classification_alter_modu}
    Fix $1\leq p,q<\infty.$ Assume that $\bigparen{\set{M_{\beta n}T_{\alpha k}g}_{k,n\in \Z},\set{M_{\beta n}T_{\alpha k}\gamma}_{k,n\in \Z}} \subseteq M^{p,q}(\R) \times (M^{p,q}(\R))^*$ is a Schauder frame for $M^{p,q}(\R).$ Then the following statements are equivalent.
    \begin{enumerate}
    \setlength\itemsep{0.5em}
        \item [\textup{(a)}] $\set{M_{\beta n}T_{\alpha k}g}_{k,n\in \Z}$ is an unconditional basis for $M^{p,q}(\R)$ plus at most finitely many elements.
        \item [\textup{(b)}] $f=\sum_{k,n\in \Z} \ip{f}{\phi_{k,n}}  M_{\beta n}T_{\alpha k}g $ with the unconditional convergence of the series in $M^{p,q}(\R)$ for all $f\in M^{p,q}(\R)$ and every alternative duals $\set{\phi_{n,k}}_{n,k\in \Z}$ of $\set{M_{\beta n}T_{\alpha k}g}_{k,n\in \Z}.$
    \end{enumerate}
\end{corollary}    
\begin{proof}
    By Theorem \ref{BS_uncon_alter_nearuc}, it suffices to show that $M^{p,q}(\R)$ does not contain a topologically isomorphic copy of $c_0.$ A special case of  \cite[Theorem 4]{BP58} states that every separable Banach space that is the dual space of some Banach space does not contain a topologically isomorphic copy of $c_0.$ The statements then follow by Lemma \ref{tf_lemmas} (f).
\end{proof}

\begin{remark}
It was shown in \cite[Theorem 3]{BCHL03} that every Gabor frame for $L^2(\R)$ is either an unconditional basis for $L^2(\R)$ or it has infinite excess. Here excess stands for the greatest number of elements that we can remove from a Gabor frame without breaking its completeness (see \cite{HY23} for the discussion on excess of general frames). Therefore, every Gabor frame is either an unconditional basis for $L^2(\R)$ or it admits at least one alternative dual for which the corresponding reconstruction formula does not converge unconditionally for all elements in $L^2(\R)$ by Corollary \ref{classification_alter_modu}. 
\end{remark}

\section*{Acknowledgments}
We thank Alexander Powell for pointing out some references and another construction of tight Gabor frames for $L^2(\R)$ with atoms not in $M^p$ for any $1\leq p<2$. We would like to express gratitude to Christopher Heil for introducing to the author the main question that was studied in Section 4 of this paper, and many thanks to him for his help and advice during the author's PhD program at Georgia Tech.

\end{document}